\documentclass[a4paper,10pt]{article}
\usepackage{amscd}
\usepackage{bbding}
\usepackage{amsfonts}
\usepackage{graphicx,amssymb,amsthm}
\usepackage{latexsym,amsmath,amscd,cite,enumerate}
\usepackage{float}
\usepackage{leftidx}
\usepackage{epstopdf}
\usepackage{manfnt}
\usepackage{stmaryrd}
\usepackage{cases}
\usepackage{bbm}
\usepackage{tabls}
\usepackage{multirow}
\usepackage{color}
\usepackage{tensor}

\allowdisplaybreaks[4]
\DeclareFontEncoding{OT2}{}{}
\DeclareTextSymbol{\cyrsftsn}{OT2}{126}
\DeclareTextSymbol{\textnumero}{OT2}{125}

\theoremstyle{definition}
\newtheorem{theorem}{Theorem}[section]
\newtheorem{lemma}{Lemma}[section]

\newtheorem{definition}{Definition}[section]
\newtheorem{remark}{Remark}[section]
\newtheorem{example}{Example}[section]

\begin{document}

\title{{The existence and stability for a class of fuzzy fractional boundary value problems constrained by variational inequalities}\thanks{This work was supported by the National Natural Science Foundation of China (11901273, 62272208, 12571318), the Natural Science Foundation of Henan (252300421997, 252300420338), the Natural Science foundation of Sichuan Province (2024NSFSC1392), the Program for Science and Technology Innovation Talents in Universities of Henan Province (23HASTIT031).}}
\author{Zengbao Wu$^{a}$\footnote{Corresponding author. E-mail addresses: zengbaowu@hotmail.com, wuzengbao@lynu.edu.cn}, Quanguo Zhang$^{a}$, Tao Chen$^{b}$, Yibin Xiao$^{c}$, Tianyin Wang$^{d}$, Chunyan Yang$^{e}$\\
$^{a}${\small \textit{Department of Mathematics, Luoyang Normal University, Luoyang, Henan  471934, P.R. China}}\\
$^{b}${\small \textit{School of Science, Southwest Petroleum University, Chengdu, Sichuan 610500, P.R. China}}\\
$^{c}${\small \textit{School of Mathematical Sciences, University of Electronic Science and Technology  of China, }}\\
{\small\textit{Chengdu, Sichuan 611731, P.R. China}}\\
$^{d}${\small \textit{School of Telecommunication Engineering, Luoyang Normal University, Luoyang, }}\\
{\small \textit{Henan  471934, P.R. China}}\\
$^{e}${\small \textit{Department of Mathematics, Sichuan University, Chengdu, Sichuan 610064, P.R. China}}
}
\date{ }
\maketitle

\begin{flushleft}
\hrulefill
\end{flushleft}
\textbf{Abstract}. In this paper, the existence and stability of a new system consisting of fuzzy fractional differential inclusions (FFDIs) with integral boundary conditions and variational inequalities are considered in Euclidean spaces. This system integrates the features of both FFDIs with integral boundary conditions and fractional differential variational inequalities in a unified framework. The existence of solutions for this system is established under some mild conditions. Moreover, the semicontinuity of the solution map for the perturbed system is investigated with respect to parameter perturbations. Finally, we use two numerical examples to show the validity of the above theoretical findings.
\begin{flushleft}
\hrulefill
\end{flushleft}

\section{Introduction}
\noindent \setcounter{equation}{0}
Let $J=[0,T]$, $K \subset R^m$ be closed, convex and nonempty, and $Q:J \times R^n \rightarrow R^m$, $S: R^m \rightarrow R^m$ be two given continuous functions. Given  $t \in J$ and $z\in R^n$, the variational inequality (VI) is to find a point $u \in K$ such that
\begin{equation}\label{VI}
  \left<Q(t,z)+S(u),v-u\right> \geq 0, \; \forall \; v \in K,
\end{equation}
where $\left<\cdot,\cdot\right>$ denotes the classical inner product in $R^m$. Let $\mbox{SOL}\left(K,Q(t,z)+S(\cdot)\right)$ denote the set of solutions to VI (\ref{VI}). This paper considers the following new system consisting of fuzzy fractional differential inclusions (FFDIs) with integral boundary conditions and VIs, which is referred to as the fuzzy fractional differential variational inequality with integral boundary conditions (abbreviated as BFFDVI):
\begin{equation}\label{BFFDVI}
\left\{
\begin{array}{l}
\tensor*[^{C}_0]{D}{_t^q} y(t) \in \left[F_{(t,y(t))}\right]_\alpha +g(t,y(t))u(t), \quad  a.e. \; t \in J,\\
u(t) \in \mbox{SOL}\left(K,Q(t,y(t))+S(\cdot)\right), \quad  a.e. \; t \in J, \\
y(0)=\int_0^T c_1(\tau,y(\tau))d\tau,\; y(T)=\int_0^T c_2(\tau,y(\tau))d\tau,
\end{array}
\right.
\end{equation}
where $q \in (1,2]$, $F:J\times R^n \rightarrow E^n$ denotes a fuzzy mapping, $g:J\times R^n \rightarrow R^{n\times m}$ and $c_i:J\times R^n \rightarrow R^n$ $(i=1,2)$ are given continuous functions. We note that when $Q(t,y(t))\equiv \mathbf{0}_{R^n}$ and $g(t,y(t))\equiv\mathbf{0}_{R^n \times R^m}$, where $\mathbf{0}_{R^n}$ and $\mathbf{0}_{R^n \times R^m}$ are the zero element in $R^n$ and $R^n \times R^m$, respectively, (\ref{BFFDVI}) degenerates into a FFDI with integral boundary conditions and a parametric VI, and the resulting FFDI with integral boundary conditions still constitutes a new class of systems. Moreover, if $q \in (0,1]$, $c_1(\tau,y(\tau))\equiv c_1$ is a constant function and $c_2(\tau,y(\tau))\equiv\mathbf{0}_{R^n}$, then (\ref{BFFDVI}) reduces to a fuzzy fractional differential variational inequality (FFDVI), which was investigated by Wu et al. \cite{WWHWW,WLZX}.

It is well known that differential variational inequalities (DVIs) belong to a family of dynamical systems, which consist of differential equations and VIs. Owing to the widespread applications in science and engineering such as microbial fermentation processes, dynamic transportation network, ideal diode circuits, dynamic Nash equilibrium problems, frictional contact problems, price control problems and so on, DVIs have currently become active areas of research. Particularly, Pang and Stewart \cite{PS} in 2008 firstly systematically investigated DVIs in Euclidean spaces. In 2015, Ke et al. \cite{K2015}, for the first time, incorporated fractional calculus into the framework of DVIs, and they investigated a fractional DVIs with delay in Euclidean spaces. In 2021, Wu et al. \cite{WWHWW} examined an FFDVI that consists of an FFDI and a VI. In 2020, Brogliato and Tanwani \cite{BT} provided an excellent review on DVIs. In 2022, Wu et al. \cite{WWHXZ} study the existence of solution and developed an approximating algorithm for a fractional differential fuzzy VI consisting of a fractional delay differential equation paired with a fuzzy VI; Zeng et al. \cite{ZCN} studied the unique solvability of a fractional differential fuzzy variational inequality with Mittag-Leffler kernel. In 2023, Mig\'{o}rski et al. \cite{MCD} examined a differential variational-hemivariational inequalities and provided an application to contact mechanics; Zhao et al. \cite{ZCL} established the existence of solutions to a differential quasi-variational-hemivariational inequality. Recently, Zeng et al. \cite{ZZH} proved the unique existence of a stochastic fractional DVI with L\'{e}vy jump and provided an application to the spatial price equilibrium problem in stochastic environments. For more works, the readers are encouraged to consult \cite{CCHX,JSZ,WCH,WQTX,zhang2023r,WCZZHX} and the citations therein.

It is worth mentioning that the fractional differential equations of order $(1,2]$ have attracted significant attention. For example, fractional Langevin equations of this order are used to characterize the super-diffusion in the anomalous diffusion of fractional Brownian motion (see, e.g., the monograph \cite{CSL}). Further details on this topic can be found in \cite{AN,BL,MK} and the references listed therein. Moreover, fractional boundary value problems, particularly those involving integral conditions, have drawn significant research interest (see, e.g.,\cite{Salem,BH2009,CH2014,ANA,AAKE}). In particular, Agarwal et al. \cite{ABH} in 2010 investigated the solvability of various classes of fractional boundary value problems. In 2015, Al-Mdallal and Hajji \cite{AH15} provided a numerical algorithm for solving nonlinear fractional boundary value problems. In 2023, Wanassi and Torres \cite{WT2023} studied a fractional differential equation with initial conditions on the function and its first derivative, along with an integral boundary condition, and they further applied the model to world population growth. Recently, Alam et al. \cite{AZA} presented the existence and stability of solutions to an implicit fractional integro-differential equation with integral boundary conditions. However, to the best of our knowledge, there are very rare works to investigate BFFDVI (\ref{BFFDVI}). The aim of our work is to make an attempt in this new direction.

The structure of the paper is as follows. Preliminaries, including notation, definitions, and lemmas, are given in Section 2. Section 3 provides a proof of the existence of solutions for BFFDVI (\ref{BFFDVI}) based on the set-valued Krasnoselskii fixed point theorem. In Section 4, we study the stability for BFFDVI (\ref{BFFDVI}). In Section 5, to validate the theoretical results, two numerical examples are presented. Section 6 concludes the paper.

\section{Preliminaries}
\noindent \setcounter{equation}{0}
This section is devoted to recalling some notation, definitions, and useful lemmas. As usual, $R^+$ denotes the set of positive reals, $L^1(J, R^n)$ denotes the totality of $R^n$-valued Lebesgue integrable functions on $J$. Denote by $L^\infty(J,R^n)$ the Banach space of bounded measurable functions $y:J\rightarrow R^n$, endowed with the norm $\|x\|_{L^\infty} = \inf \{c>0: \|y(t)\| \leq c, \; a.e. \; t\in J\}$, and $C(J, R^n)$ denotes the totality of $R^n$-valued continuous functions on $J$ with the norm $\|y\|_{C}=\underset{t\in J}\max\|y(t)\|$.

Let $Y$ be a base space, $\mathcal{F}(Y)$ be the set of all fuzzy sets of $Y$. A fuzzy set is defined as a mapping $\omega : Y \rightarrow [0,1]$. A fuzzy mapping is a mapping $F:\Lambda \rightarrow \mathcal{F}(Y)$, where $\emptyset\neq\Lambda \subset Y$. For each $y \in \Lambda$, the image $F(y)$, denoted by $F_y$, is a fuzzy set, and $F_y(\theta)$ denotes membership grade of $\theta$ in $F_y$. For $\alpha \in (0,1]$, the $\alpha$-level set of $\omega$ is $[\omega]_\alpha = \{\theta \in Y: \omega(\theta)\geq\alpha\}$. The support of $\omega$ is defined as $[\omega]_0=\overline{\underset{\alpha \in (0,1]}\cup [\omega]_\alpha}$, where the overline denotes the closure of the union. Let us denote by $E^n$ the space consisting of all fuzzy sets of $R^n$ satisfying normal, fuzzy convex, upper semicontinuous as function with compact level sets (see, e.g., \cite[p.5]{PD}).

Given two Banach spaces $Z_1, Z_2$. A set-valued mapping $\Upsilon: Z_1 \rightarrow 2^{Z_2}\setminus \{\emptyset\}$ has convex (compact, closed) values if $\Upsilon(z)$ is convex (compact, closed) for all $z \in Z_1$. $\Upsilon$ has a fixed point if there exists a point $z \in Z_1 \subset Z_2$ such that $z \in \Upsilon(z)$. $\Upsilon$ is called upper semicontinuous (u.s.c.) on $Z_1$ provided that for any $z_0 \in Z_1$ and every open set $U \subset Z_2$ containing $\Upsilon(z_0)$, one can find an open neighbourhood $O$ of $z_0$ such that $\Upsilon(z) \subset U$ for all $z \in O$. $\Upsilon$ is called lower semicontinuous (l.s.c.) on $Z_1$ whenever for each $z_0 \in Z_1$ and any open set $U \subset Z_2$ satisfying $\Upsilon(z_0) \cap U \neq \emptyset$, there exists an open neighbourhood $O$ of $z_0$ ensuring $\Upsilon(z) \cap U \neq \emptyset$ for every $z \in O$. We say that $\Upsilon$ is continuous if it satisfies both upper and lower semicontinuity. $\Upsilon$ is called completely continuous whenever $\Upsilon(U)$ is relatively compact for each bounded subset $U \subset Z_1$. $\Upsilon$ has a closed graph if the graph $Gr(\Upsilon) = \{(z_1, z_2) \in Z_1 \times Z_2 : z_2 \in \Upsilon(z_1)\}$ of $\Upsilon$ is a closed set of $Z_1 \times Z_2$. It is noted that for any completely continuous mapping $\Upsilon$ with compact values, upper semicontinuity is equivalent to the mapping possessing a closed graph (see \cite{BHNO}).
\begin{definition}\cite{KD}
The Riemann-Liouville fractional integral is defined as
  \begin{equation*}
\tensor*{I}{_0^q}y(t) = \frac{1}{\Gamma(q)}\int_0^t(t-\tau)^{q-1}y(\tau)d\tau, \; q>0
\end{equation*}
where the Gamma function $\Gamma$ is defined by $\displaystyle \Gamma(q)=\int_{0}^\infty \tau^{q-1}e^{-\tau}d\tau$. Notably, $\Gamma(q+1)=q\Gamma(q)$ and $\Gamma(1)=1$.
\end{definition}

\begin{definition}\cite{KD}
The Caputo fractional order derivative is defined as
      \begin{eqnarray*}
\tensor*[^{C}_0]{D}{_t^q}y(t) = \frac{1}{\Gamma(n-q)}\int_0^t(t-\tau)^{n-q-1}y^{(n)}(\tau)d \tau,
\end{eqnarray*}
here $n=[q]+1$ with $[q]$ denoting the floor of $q$.
\end{definition}

Motivated by \cite[Definition 2.4]{WWHWW} and \cite[Lemma 3.21]{ABH}, we introduce the definition of BFFDVI (\ref{BFFDVI}) as follows.

\begin{definition}\label{def-BFFDVI} For $y\in C(J,R^n)$ and a integrable function $u:J \rightarrow K$, we say that $(y,u)$ is a mild solution of BFFDVI (\ref{BFFDVI}) if
\begin{equation*}
\left\{
\begin{array}{l}
\displaystyle y(t) = \frac{1}{\Gamma(q)} \int_0^t(t-\tau)^{q-1}\left[f(\tau) + g(\tau,y(\tau))u(\tau)\right]d\tau- \frac{t}{\Gamma(q)T}\int_0^T(T-\tau)^{q-1}\left[f(\tau) + g(\tau,y(\tau))u(\tau)\right]d\tau\\
\displaystyle\qquad \quad +\frac{t}{T}\int_0^T c_2(\tau,y(\tau))d\tau + \left(1-\frac{t}{T}\right)\int_0^T c_1(\tau,y(\tau))d\tau, \quad t \in J,\\
u(t) \in \mbox{SOL}\left(K,Q(t,y(t))+S(\cdot)\right), \quad a.e. \; t \in J,
\end{array}
\right.
\end{equation*}
where $f \in \mathbb{S}^1_{\widetilde{F}}(y)$ and
\begin{equation}\label{S1F}
\mathbb{S}^1_{\widetilde{F}}(y)=\left\{z \in L^1(J,R^n): z(\tau) \in  \widetilde{F}(\tau,y(\tau)) = \left[F_{(\tau,y(\tau))}\right]_\alpha, \; a.e. \; \tau\in J\right\}.
\end{equation}
Here, $u$ and $y$ are termed the variational control trajectory and the mild trajectory, respectively.
\end{definition}

\begin{remark}\label{Def-Rem}
Let
\begin{equation}\label{GQS}
G(t,y(t))=\left\{u(t): u(t) \in \mbox{SOL}(K,Q(t,y(t))+S(\cdot))\right\}.
\end{equation}
It follows from Definition \ref{def-BFFDVI} that the existence of mild solution of GFFDVI (\ref{BFFDVI}) can be reformulated by the existence of the following system
\begin{align*}
y(t)=&\frac{1}{\Gamma(q)} \int_0^t(t-\tau)^{q-1}\left[f(\tau) + g(\tau,y(\tau))h(\tau)\right]d\tau- \frac{t}{\Gamma(q)T}\int_0^T(T-\tau)^{q-1}\left[f(\tau) + g(\tau,y(\tau))h(\tau)\right]d\tau\\
& +\frac{t}{T}\int_0^T c_2(\tau,y(\tau))d\tau + \left(1-\frac{t}{T}\right)\int_0^T c_1(\tau,y(\tau))d\tau, \quad t \in J,
\end{align*}
where $f \in \mathbb{S}^1_{\widetilde{F}}(y)$ and $h \in \mathbb{S}^1_{G}(y)$ with $\mathbb{S}^1_{\widetilde{F}}(y)$ being defined by (\ref{S1F}) and
\begin{equation}\label{S1G}
\mathbb{S}^1_{G}(y)=\left\{z \in L^1(J,R^n): z(\tau) \in  G(\tau,y(\tau)), \; a.e. \; \tau \in J \right\}.
\end{equation}
\end{remark}

The following lemma is obtained as an immediate consequence of \cite[Lemma 8.6.4]{Chang} and \cite[Lemma 5.3.5(iii)]{Chang}.

\begin{lemma}\label{Bothmea}
Consider two measurable set-valued mappings $\Lambda_{1}, \Lambda_{2}:J\rightarrow 2^{R^n}$ with compact values. Suppose $g_1:J\rightarrow R^n$ is a measurable selection of $\Lambda_1$. Then there exists another measurable selection $g_2:J\rightarrow R^n$ of $\Lambda_2$ satisfying
$$\|g_2(t)-g_1(t)\|\leq H(\Lambda_2(t),\Lambda_1(t))$$
for all $t\in J$ with $H$ denoting the Hausdorff distance.
\end{lemma}

\begin{lemma}\cite[Theorem 5]{Nadler} \label{Fixpoint}
Assume that the metric space $(Z,d)$ is complete. Consider a set-valued mapping $\Phi:Z \rightarrow 2^Z$ whose values are closed and bounded. Suppose that $\Phi$ is contractive in the sense of the Hausdorff distance $H$: there is $\varrho \in (0,1)$ such that
$$H(\Phi(\zeta_1),\Phi(\zeta_2))\leq \varrho d(\zeta_1,\zeta_2), \quad \forall~\zeta_1, \zeta_2 \in Z.$$
Then $\Phi$ admits a fixed point in $Z$.
\end{lemma}

\begin{lemma}\cite[Corollary 3.2]{Dhage}\label{K-fix}
(Set-valued Krasnoselskii fixed point theorem) Let $\overline{B_\delta(0)}$ represent the closed ball of radius $\delta>0$ centered at the origin in a Banach space $X$. Consider two set-valued mappings
$$A:\overline{B_\delta(0)} \rightarrow 2^X\setminus \{\emptyset\}, \;\; B:\overline{B_\delta(0)} \rightarrow 2^X\setminus \{\emptyset\},$$
where $A$ has convex, bounded and closed values, while $B$ has convex and compact values. Assume further that $A$ is contractive with respect to the Hausdorff distance, and that $B$ is completely continuous and u.s.c. Then the following alternative holds: either
 \begin{enumerate}[(i)]
  \item the sum $A +B$ admits a fixed point in $\overline{B_r(0)}$, or
  \item there exists $z \in X$ with $\|z\|=\delta$ and a scalar $\kappa>1$ such that $\kappa z \in Az +Bz$.
\end{enumerate}
\end{lemma}

\section{Existence results for BFFDVI}
\noindent \setcounter{equation}{0}
This section is devoted to the existence of solutions of BFFDVI (\ref{BFFDVI}). In the sequel, we assume that:
\begin{enumerate}
  \item[(A$_1$)] $\mathcal{H}\left(F_{(t,z_1)},F_{(t,z_2)}\right)\leq L_F\|z_1-z_2\| \; (L_F>0)$  $\forall \; t \in J, \; z_1, z_2 \in R^n$, where $\mathcal{H}$ is a metric on $E^n$ (see, e.g., \cite{LM}) and is defined by
$$\mathcal{H}(\omega_1,\omega_2)=\sup\left\{H([\omega_1]_\alpha,[\omega_2]_\alpha): 0 \leq \alpha \leq 1\right\}\quad \mbox{for all} \; \omega_1,\omega_2 \in E^n$$
with $H$ being the Hausdorff distance between two sets;
  \item[(A$_2$)] for every $z \in R^n$, $F_{(\cdot,z)}$ is strongly measurable;
  \item[(A$_3$)] for every $z\in R^n$ and $a.e.$ $t \in J$, it holds $\|F_{(t,z)}\| \leq p(t)$, where $p \in L^\infty(J,R^+)$;
  \item[(A$_4$)] there exists $\eta_g>0$ satisfying $\|g(t,z)\|\leq \eta_g$ for all $t\in J$, $z\in R^n$;
  \item[(A$_5$)] there exists $\eta_Q>0$ satisfying $\|Q(t,z)\|\leq \eta_Q$ for all $t\in J$, $z\in R^n$;
  \item[(A$_6$)] there is $u_0 \in K$ satisfying
  $$\liminf \limits_{u \in K, \|u\| \rightarrow \infty} \frac{\left<S(u),u-u_0\right>}{\|u\|^2} >0.$$
  Moreover, $S$ is monotone on $K$.
\end{enumerate}

\begin{remark}\label{R-assume}
Assumption (A$_1$) is weaker than \cite[Hypothesis (H$_1$)]{WWHWW}. Assumption (A$_3$) is a special case of \cite[Hypothesis (H$_3$)]{WWHWW}. Assumptions (A$_2$) and (A$_4$)-(A$_6$) are the same as\cite[Hypotheses (H$_2$) and (H$_4$)-(H$_6$)]{WWHWW}.
\end{remark}

\begin{lemma}\label{F-usc}
Let $F:J \times R^n \rightarrow E^n$ be a fuzzy mapping and $\widetilde{F}:J\times R^n \rightarrow 2^{R^n}$ be defined by
\begin{equation}\label{F}
\widetilde{F}(t,z) = \left[F_{(t,z)}\right]_\alpha=\left\{x\in R^n: F_{(t,z)}(x) \geq \alpha\right\},
\end{equation}
where $\alpha \in [0,1]$, $z \in R^n$. Then $\widetilde{F}$ has nonempty convex and compact values. Furthermore, if assumption (A$_1$) holds, then, for any $t \in J$, $\widetilde{F}(t,\cdot)$ is Lipschitz continuous in its second argument. In addition, for every $x \in \widetilde{F}(t,z)$, we have
\begin{equation}\label{xFy}
\|x\| \leq L_F \|z\|+\left\|\widetilde{F}(t,0)\right\|, \quad \forall \; t \in J, \; z \in R^n.
\end{equation}

\textbf{Proof.}\hspace{0.2cm} Combining \cite[Lemma 3.1]{WWHWW} and \cite[Lemma 3.1]{WLZX}, it follows immediately that Lemma \ref{F-usc} holds.
\end{lemma}

\begin{remark}\label{Fpr}
In light of assumption (A$_3$) and (3.3) in \cite[Lemma 3.4]{WWHWW}, for any $z \in R^n$ and $a.e.$ $t \in J$, one has
\begin{equation*}
\sup\left\{\|x\|: x \in \widetilde{F}(t,z) = \left[F_{(t,z)}\right]_\alpha\right\} \leq  p(t).
\end{equation*}
\end{remark}

Similar to the argument of \cite[Lemma 3.4]{WWHWW}, we conclude that the following result holds.

\begin{lemma}\label{FG-nonempty}
Let (A$_1$)-(A$_6$) hold. Then $\mathbb{S}^1_{\widetilde{F}}(y) \neq \emptyset$ and $\mathbb{S}^1_{G}(y) \neq \emptyset$, where $\mathbb{S}^1_{\widetilde{F}}(y), \mathbb{S}^1_{G}(y)$ are defined by (\ref{S1F}) and (\ref{S1G}), respectively.
\end{lemma}

\begin{remark}\label{G-usc}
Let (A$_5$)-(A$_6$) hold. The set-valued mapping $G:J \times R^n \rightarrow 2^{R^n}$ given by (\ref{GQS}) is u.s.c., and each of its images is nonempty convex and compact. In addition, for every $t\in J$, $z\in R^n$,
$$\|G(t,z)\|=\sup\{\|x\|:x\in G(t,z)\} \leq \eta_S(1+\|Q(t,z)\|) \leq  \eta_S(1+\eta_Q),$$
where $\eta_S >0$ is a constant (see \cite[Remark 3.2]{WWHWW}).
\end{remark}

\begin{lemma}\label{phi-Phi}
Let (A$_1$)-(A$_3$) hold. Then the function $\varphi:J \rightarrow R^n$ by setting
$$
\varphi(t) = \frac{1}{\Gamma(q)} \int_0^t(t-\tau)^{q-1}f(\tau) d\tau- \frac{t}{\Gamma(q)T}\int_0^T(T-\tau)^{q-1}f(\tau)d\tau
$$
is continuous, where $f \in \mathbb{S}^1_{\widetilde{F}}(y)$, $y \in C(J,R^n)$. Moreover, the set-valued mapping $\Phi:C(J,R^n) \rightarrow 2^{C(J,R^n)}$ given by
\begin{equation}\label{Phi}
 \Phi(y) = \left\{\varphi \in C(J,R^n) : \varphi(t)= \frac{1}{\Gamma(q)} \int_0^t(t-\tau)^{q-1}f(\tau) d\tau- \frac{t}{\Gamma(q)T}\int_0^T(T-\tau)^{q-1}f(\tau)d\tau, \; f \in \mathbb{S}^1_{\widetilde{F}}(y)\right\}
\end{equation}
has bounded, closed and convex values. Furthermore, if $\displaystyle\rho=\frac{2 L_F T^q}{\Gamma(q+1)} <1$, then $\Phi$ is contractive.

\textbf{Proof.}\hspace{0.2cm} The proof proceeds in three steps.

{\bf Step 1.} We establish the continuity of $\varphi$.

By virtue of Lemma \ref{FG-nonempty}, $\varphi$ is well defined.  Let $0\leq t_1 < t_2 \leq T$. Given $y \in C(J,R^n)$, $f \in \mathbb{S}^1_{\widetilde{F}}(y)$, we have
\begin{align}\label{phi-1}
\varphi(t_2) - \varphi(t_1)
 =& \frac{1}{\Gamma(q)} \int_{t_1}^{t_2}(t_2-\tau)^{q-1}f(\tau) d\tau +\frac{1}{\Gamma(q)} \int_0^{t_1}\left[(t_2-\tau)^{q-1}-(t_1-\tau)^{q-1}\right]f(\tau) d\tau\nonumber \\
& - \frac{t_2-t_1}{\Gamma(q)T}\int_0^T(T-\tau)^{q-1}f(\tau)d\tau.
\end{align}
Combining assumption (A$_3$), Remark \ref{Fpr} and (\ref{phi-1}), one has
\begin{eqnarray*}
 && \|\varphi(t_2) - \varphi(t_1)\| \nonumber \\
 &\leq&  \frac{1}{\Gamma(q)} \int_{t_1}^{t_2}(t_2-\tau)^{q-1}\|f(\tau)\| d\tau +\frac{1}{\Gamma(q)} \int_0^{t_1}\left|(t_2-\tau)^{q-1}-(t_1-\tau)^{q-1}\right|\|f(\tau)\| d\tau \nonumber \\
  && + \frac{t_2-t_1}{\Gamma(q)T}\int_0^T(T-\tau)^{q-1}\|f(\tau)\|d\tau\nonumber \\
   &\leq&   \frac{\|p\|_{L^\infty}}{\Gamma(q)} \int_{t_1}^{t_2}(t_2-\tau)^{q-1} d\tau +\frac{\|p\|_{L^\infty}}{\Gamma(q)} \int_0^{t_1}\left[(t_2-\tau)^{q-1}-(t_1-\tau)^{q-1}\right] d\tau\nonumber \\
  &&+ \frac{(t_2-t_1)\|p\|_{L^\infty}}{\Gamma(q)T} \int_0^T(T-\tau)^{q-1}d\tau \nonumber \\
 &=&\frac{\|p\|_{L^\infty}T^{q-1}}{\Gamma(q+1)} (t_2-t_1) + \frac{\|p\|_{L^\infty}}{\Gamma(q+1)} \left(t^q_2-t^q_1\right).
\end{eqnarray*}
Thus $\varphi$ is continuous.

{\bf Step 2.} We show that $\Phi(y)$ is a bounded, convex and closed set for any given $y \in C(J,R^n)$.

For any $\varphi \in \Phi(y)$, one has
$$
\varphi(t) = \frac{1}{\Gamma(q)} \int_0^t(t-\tau)^{q-1}f(\tau) d\tau- \frac{t}{\Gamma(q)T}\int_0^T(T-\tau)^{q-1}f(\tau)d\tau
$$
with $f \in \mathbb{S}^1_{\widetilde{F}}(y)$. Combining assumption (A$_3$) with Remark \ref{Fpr} yields
\begin{eqnarray*}
 \|\varphi(t) \|
  &\leq& \frac{1}{\Gamma(q)} \int_0^t(t-\tau)^{q-1}\|f(\tau)\| d\tau+ \frac{1}{\Gamma(q)}\int_0^T(T-\tau)^{q-1}\|f(\tau)\|d\tau \nonumber \\
   &\leq& \frac{\|p\|_{L^\infty}}{\Gamma(q)} \int_0^t(t-\tau)^{q-1} d\tau + \frac{\|p\|_{L^\infty}}{\Gamma(q)}\int_0^T(T-\tau)^{q-1}d\tau \nonumber \\
 &\leq& 2\frac{\|p\|_{L^\infty}}{\Gamma(q+1)} T^q.
\end{eqnarray*}
As a result, $\Phi(y)$ remains bounded.

Next, we show that $\Phi(y)$ is a convex set. Let $\varphi_1,\varphi_2 \in \Phi(y)$. Then there exists $f_1,f_2 \in \mathbb{S}^1_{\widetilde{F}}(y)$ such that
$$
\varphi_i(t) = \frac{1}{\Gamma(q)} \int_0^t(t-\tau)^{q-1}f_i(\tau) d\tau- \frac{t}{\Gamma(q)T}\int_0^T(T-\tau)^{q-1}f_i(\tau)d\tau \;\; (i=1,2).
$$
It follows that, for any $0\leq \varrho \leq 1$,
\begin{eqnarray*}
 && \varrho\varphi_1(t)+(1-\varrho)\varphi_2(t) \nonumber \\
  &=&  \frac{1}{\Gamma(q)} \int_0^t(t-\tau)^{q-1}\left[\varrho f_1(\tau)+(1-\varrho)f_2(\tau)\right] d\tau- \frac{t}{\Gamma(q)T}\int_0^T(T-\tau)^{q-1}\left[\varrho f_1(\tau)+(1-\varrho)f_2(\tau)\right]d\tau.
\end{eqnarray*}
The convex-valuedness of $\widetilde{F}$ implies that of  $\mathbb{S}^1_{\widetilde{F}}(y)$ (see Remark 2.1 in \cite{DGO}). Thus $\varrho f_1+(1-\varrho)f_2 \in \mathbb{S}^1_{\widetilde{F}}(y)$. Consequently, $\varrho \varphi_1+ (1-\varrho)\varphi_2 \in \Phi(y)$, that is, $\Phi(y)$ is a convex set.

We finish by showing the closedness of $\Phi(y)$. Let $\{\varphi_n\} \subset \Phi(y)$ be a sequence with $\varphi_n \rightarrow \varphi$. We have
\begin{equation}\label{phin}
\varphi_n(t) = \frac{1}{\Gamma(q)} \int_0^t(t-\tau)^{q-1}f_n(\tau) d\tau- \frac{t}{\Gamma(q)T}\int_0^T(T-\tau)^{q-1}f_n(\tau)d\tau,
\end{equation}
where $\{f_n\} \subset \mathbb{S}^1_{\widetilde{F}}(y)$, $n=1,2,\ldots.$ From Remark \ref{Fpr}, $\|f_n(\tau)\| \leq p(\tau)$ $a.e.$ on $J$. Thus $\{f_n: n\geq1\}$ is integrably bounded. By \cite[Corollary 13, Section 19.5]{RF}, some subsequence, again denoted $\{f_n\}$, converges weakly to $\widetilde{f} \in L^1([0,T],R^n)$. Mazur's Lemma (see Lemma A.3 in \cite{Stewart}) then yields convex combinations 
\begin{equation}\label{x-rho-f}
x_n=\sum\limits_{j=n}^{j_0(n)}\rho_{n,j}f_j \rightarrow \widetilde{f} \in L^1([0,T],R^n),  
\end{equation}
with $j_0(n) \in \mathbb{N}$, $j_0(n)>n$, $\rho_{n,j}\geq0$, and $\sum\limits_{j=n}^{j_0(n)}\rho_{n,j}=1$. As $x_n \rightarrow \widetilde{f}$ in $L^1([0,T],R^n)$, we may assume $x_n(\tau) \rightarrow \widetilde{f}(\tau)$ $a.e.$ on $J$. Lemma \ref{F-usc} ensures that $\widetilde{F}$ has the convex and closed values. Thus, for $a.e.$ $\tau \in J$
\begin{equation*}
\widetilde{f}(\tau) \in \mathop{\bigcap}\limits_{n \geq1} \overline{\{x_j(\tau):j\geq n\}} \subseteq \mathop{\bigcap}\limits_{n \geq1} \overline{\mbox{conv}}\{f_j(\tau): j\geq n\} \subseteq \widetilde{F}(\tau,y(\tau)),
\end{equation*}
where the overline denotes the closure. Hence $f \in \mathbb{S}^1_{\widetilde{F}}(y)$.
Let $\widetilde{\varphi}_n=\sum\limits_{j=n}^{j_0(n)}\rho_{n,j}\varphi_j$. Then $\widetilde{\varphi}_n(t) \rightarrow \varphi(t)$ for every $t\in J$. In view of (\ref{phin}), one has
\begin{equation}\label{phin1}
\widetilde{\varphi}_n(t) = \frac{1}{\Gamma(q)} \int_0^t(t-\tau)^{q-1}x_n(\tau) d\tau- \frac{t}{\Gamma(q)T}\int_0^T(T-\tau)^{q-1}x_n(\tau)d\tau.
\end{equation}
It is note that for every $t \in J$, $\tau \in (0,t]$,
$$\|x_n(\tau)\| \leq  p(\tau) \quad \mbox{and} \quad \left\|(t-\tau)^{q-1}x_n(\tau)\right\| \leq \left(t-\tau\right)^{q-1} p(\tau).$$
It follows from assumption (A$_3$) that $p(\cdot)\in L^1(J,R^+)$ and $\left(t-\cdot\right)^{q-1} p(\cdot) \in L^1(J,R^+)$. By passing to the limit as $n \rightarrow \infty$ in (\ref{phin1}), we have
$$
\varphi(t) = \frac{1}{\Gamma(q)} \int_0^t(t-\tau)^{q-1}f(\tau) d\tau- \frac{t}{\Gamma(q)T}\int_0^T(T-\tau)^{q-1}f(\tau)d\tau,
$$
where $f \in \mathbb{S}^1_{\widetilde{F}}(y)$. Accordingly, the values of $\Phi$ are closed.

{\bf Step 3.} We claim that $\Phi$ satisfies the contractive property.

For any $y_1,y_2 \in C(J,R^n)$ and $\varphi_1 \in \Phi(y_1)$, there is $f_1 \in \mathbb{S}^1_{\widetilde{F}}(y_1)$ such that
\begin{equation}\label{varphi1}
\varphi_1(t) = \frac{1}{\Gamma(q)} \int_0^t(t-\tau)^{q-1}f_1(\tau) d\tau- \frac{t}{\Gamma(q)T}\int_0^T(T-\tau)^{q-1}f_1(\tau)d\tau.
\end{equation}
It follows from assumptions (A$_1$)-(A$_2$) and Lemma \ref{F-usc} that $\widetilde{F}$ is compact-valued and satisfies the Carath\'{e}odory conditions (cf. \cite[Definition 1.3.5]{KOZ}). By virtue of \cite[Theorem 1.3.4]{KOZ}, we obtain the measurability of $\widetilde{F}(\cdot,y_1(\cdot))$ and $\widetilde{F}(\cdot,y_2(\cdot))$. Lemma \ref{Bothmea} then allows us to select a measurable $f_2(\tau) \in \widetilde{F}(\tau,y_2(\tau))$ such that
\begin{eqnarray}\label{f12F}
&&\|f_1(\tau)-f_2(\tau)\| \nonumber\\
&\leq& H\left(\widetilde{F}(\tau,y_1(\tau)),\widetilde{F}(\tau,y_2(\tau))\right)
=H\left(\left[F_{(t,y_1)}\right]_\alpha, \left[F_{(t,y_2)}\right]_\alpha\right) \nonumber\\
&\leq& \mathcal{H}\left(F_{(t,y_1)},F_{(t,y_2)}\right)
\leq L_F \|y_1(\tau)-y_2(\tau)\|
\end{eqnarray}
for all $\tau \in J$. In view of Remark \ref{Fpr}, one has $f_2 \in L^\infty$ and so $f_2$ is Lebesgue integrable on $J$. Consequently, $f_2 \in \mathbb{S}^1_{\widetilde{F}}(y_2)$. Let
\begin{equation}\label{varphi2}
\varphi_2(t) = \frac{1}{\Gamma(q)} \int_0^t(t-\tau)^{q-1}f_2(\tau) d\tau- \frac{t}{\Gamma(q)T}\int_0^T(T-\tau)^{q-1}f_2(\tau)d\tau.
\end{equation}
Then $\varphi_2 \in \Phi(y_2)$. Combining (\ref{varphi1}), (\ref{f12F}) and (\ref{varphi2}), we have
\begin{eqnarray*}
&& \left\|\varphi_1(t) - \varphi_2(t)\right\|   \nonumber\\
&\leq& \frac{1}{\Gamma(q)} \int_0^t(t-\tau)^{q-1}\|f_1(\tau)-f_2(\tau)\| d\tau + \frac{1}{\Gamma(q)} \int_0^T(T-\tau)^{q-1}\|f_1(\tau)-f_2(\tau)\| d\tau \nonumber\\
&\leq& \frac{L_F}{\Gamma(q)} \int_0^t(t-\tau)^{q-1}\|y_1(\tau)-y_2(\tau)\|d\tau  + \frac{L_F}{\Gamma(q)} \int_0^T(T-\tau)^{q-1}\|y_1(\tau)-y_2(\tau)\|d\tau \nonumber\\
&\leq&\frac{L_F}{\Gamma(q)} \|y_1-y_2\|_C \int_0^t(t-\tau)^{q-1}d\tau  + \frac{L_F}{\Gamma(q)}\|y_1-y_2\|_C  \int_0^T(T-\tau)^{q-1}d\tau \nonumber\\
&\leq&\frac{2 L_F T^q}{\Gamma(q+1)}\|y_1-y_2\|_C  \nonumber\\
&=& \rho \left\|y_1-y_2\right\|_C.
\end{eqnarray*}
Consequently,
$$
\left\|\varphi_1 - \varphi_2\right\|_C
\leq \rho \left\|y_1-y_2\right\|_C.
$$
It yields
\begin{equation*}
d(\varphi_1,\Phi(y_2)) = \underset{\varphi_2 \in \Phi(y_2)}\inf\| \varphi_1-\varphi_2\|_{C}
\leq \rho \left\|y_1-y_2\right\|_{C}.
\end{equation*}
Since $\varphi_1 \in \Phi(y_1)$ is arbitrary, one gets
\begin{equation*}
\underset{\varphi_1 \in \Phi(y_1)}\sup
d(\varphi_1,\Phi(y_2)) \leq \rho \left\|y_1-y_2\right\|_{C}.
\end{equation*}
Similarly, we have
\begin{equation*}
\underset{\varphi_2 \in \Phi(y_2)}\sup d(\varphi_2,\Phi(y_1)) \leq \rho \left\|y_1-y_2\right\|_{C}.
\end{equation*}
It follows that
\begin{equation*}
H\left(\Phi(y_1),\Phi(y_2)\right) \leq \rho \left\|y_1-y_2\right\|_{C},
\end{equation*}
which yields that $\Phi$ is contractive since $\rho<1$. \hfill $\Box$
\end{lemma}

\begin{lemma}\label{psi-Psi}
Let (A$_4$)-(A$_6$) hold. If there exist $M_{1}, M_{2}>0$ such that $\|c_1(t,z)\|\leq M_{1}$ and $\|c_2(t,z)\|\leq M_{2}$ for any $t\in J$, $z\in R^n$, then the function $\psi:J \rightarrow R^n$ by setting
\begin{align*}
\psi(t) =& \frac{1}{\Gamma(q)} \int_0^t(t-\tau)^{q-1}g(\tau,y(\tau))h(\tau)d\tau- \frac{t}{\Gamma(q)T}\int_0^T(T-\tau)^{q-1}g(\tau,y(\tau))h(\tau)d\tau \nonumber \\
&+\frac{t}{T}\int_0^T  c_2(\tau,y(\tau))d\tau + \left(1-\frac{t}{T}\right)\int_0^T c_1(\tau,y(\tau))d\tau
\end{align*}
is continuous, where $h \in \mathbb{S}^1_G(y)$, $y \in C(J,R^n)$. Moreover, the set-valued mapping $\Psi:C(J,R^n) \rightarrow 2^{C(J,R^n)}$ given by
\begin{align}\label{Psi}
 \Psi(y)=&  \left\{\psi \in C(J,R^n) : \psi(t)= \frac{t}{T}\int_0^T  c_2(\tau,y(\tau))d\tau + \left(1-\frac{t}{T}\right)\int_0^T c_1(\tau,y(\tau))d\tau \right.  \nonumber \\
   & \left. +\frac{1}{\Gamma(q)} \int_0^t(t-\tau)^{q-1}g(\tau,y(\tau))h(\tau) d\tau- \frac{t}{\Gamma(q)T}\int_0^T(T-\tau)^{q-1}g(\tau,y(\tau))h(\tau)d\tau, \; h \in \mathbb{S}^1_G(y)\right\}
\end{align}
has compact and convex values. In addition, $\Psi$ is completely continuous and is u.s.c.

\textbf{Proof. }\hspace{0.2cm} The proof proceeds in three steps.

{\bf Step 1.} We show the continuity of $\psi$.

By Lemma \ref{FG-nonempty}, we have that $\psi$ is well defined.  Let $0\leq t_1 < t_2 \leq T$. Given $y \in C(J,R^n)$, $h \in \mathbb{S}^1_G(y)$, one has
\begin{eqnarray*}
 && \psi(t_2) - \psi(t_1) \nonumber \\
 &=& \frac{1}{\Gamma(q)} \int_{t_1}^{t_2}(t_2-\tau)^{q-1}g(\tau,y(\tau))h(\tau) d\tau +\frac{1}{\Gamma(q)} \int_0^{t_1}\left[(t_2-\tau)^{q-1}-(t_1-\tau)^{q-1}\right]g(\tau,y(\tau))h(\tau) d\tau\nonumber \\
&& - \frac{t_2-t_1}{\Gamma(q)T}\int_0^T(T-\tau)^{q-1}g(\tau,y(\tau))h(\tau)d\tau + \frac{t_2-t_1}{T}\int_0^T  \left[c_2(\tau,y(\tau))-c_1(\tau,y(\tau)) \right]d\tau.
\end{eqnarray*}
It follows from assumption (A$_4$) and Remark \ref{G-usc} that
\begin{eqnarray}\label{psi-2}
 && \|\psi(t_2) - \psi(t_1)\| \nonumber \\
 &\leq&  \frac{1}{\Gamma(q)} \int_{t_1}^{t_2}(t_2-\tau)^{q-1}\|g(\tau,y(\tau))\|\|h(\tau)\|d\tau +\frac{1}{\Gamma(q)} \int_0^{t_1}\left|(t_2-\tau)^{q-1} -(t_1-\tau)^{q-1}\right|\|g(\tau,y(\tau))\|\|h(\tau)\|d\tau \nonumber \\
  && + \frac{t_2-t_1}{\Gamma(q)T}\int_0^T(T-\tau)^{q-1} \|g(\tau,y(\tau))\|\|h(\tau)\|d\tau + \frac{t_2-t_1}{T}\int_0^T  \left\|c_2(\tau,y(\tau))-c_1(\tau,y(\tau)) \right\|d\tau\nonumber \\
   &\leq&   \frac{\eta_g \eta_S(1+\eta_Q)}{\Gamma(q)} \int_{t_1}^{t_2}(t_2-\tau)^{q-1} d\tau +\frac{\eta_g \eta_S(1+\eta_Q)}{\Gamma(q)} \int_0^{t_1}\left[(t_2-\tau)^{q-1}-(t_1-\tau)^{q-1}\right] d\tau\nonumber \\
  &&+ \frac{(t_2-t_1)\eta_g \eta_S(1+\eta_Q)}{\Gamma(q)T} \int_0^T(T-\tau)^{q-1}d\tau + (M_1+M_2)(t_2-t_1) \nonumber \\
 &=&\left(\frac{\eta_g \eta_S(1+\eta_Q)T^{q-1}}{\Gamma(q+1)}+M_1+M_2\right) (t_2-t_1) + \frac{\eta_g \eta_S(1+\eta_Q)}{\Gamma(q+1)} \left(t^q_2-t^q_1\right).
\end{eqnarray}
Hence $\psi$ is continuous.

{\bf Step 2.} We show that $\Psi$ possesses the property of complete continuity, and its values are both compact and convex.

Let $\Omega \subset C(J, R^n)$ be a bounded set. We claim that $\Psi(\Omega)$ is a uniformly bounded. In fact, for any $\psi \in \Psi(\Omega)=\underset{y\in \Omega}\bigcup\Psi(y)$, there is $y \in \Omega$ and $h \in \mathbb{S}^1_G(y)$ such that
\begin{align*}
\psi(t) =& \frac{1}{\Gamma(q)} \int_0^t(t-\tau)^{q-1}g(\tau,y(\tau))h(\tau) d\tau- \frac{t}{\Gamma(q)T}\int_0^T(T-\tau)^{q-1}g(\tau,y(\tau))h(\tau)d\tau \nonumber \\
&+\frac{t}{T}\int_0^T  c_2(\tau,y(\tau))d\tau + \left(1-\frac{t}{T}\right)\int_0^T c_1(\tau,y(\tau))d\tau, \quad t \in J.
\end{align*}
Using assumption (A$_4$) and Remark \ref{G-usc}, one has
\begin{eqnarray*}
 \|\psi(t) \|
  &\leq& \frac{1}{\Gamma(q)} \int_0^t(t-\tau)^{q-1}\|g(\tau,y(\tau))\|\|h(\tau)\|d\tau+ \frac{1}{\Gamma(q)}\int_0^T(T-\tau)^{q-1}\|g(\tau,y(\tau))\|\|h(\tau)\|d\tau \nonumber \\
&&+ \int_0^T  \|c_2(\tau,y(\tau))\|d\tau + \int_0^T \|c_1(\tau,y(\tau))\|d\tau \nonumber \\
   &\leq& \frac{\eta_g \eta_S(1+\eta_Q)}{\Gamma(q)} \int_0^t(t-\tau)^{q-1} d\tau + \frac{\eta_g \eta_S(1+\eta_Q)}{\Gamma(q)}\int_0^T(T-\tau)^{q-1}d\tau+(M_1+M_2)T \nonumber \\
 &\leq& 2\frac{\eta_g \eta_S(1+\eta_Q)}{\Gamma(q+1)} T^q+(M_1+M_2)T.
\end{eqnarray*}
This yields the uniform boundedness of $\Psi(\Omega)$. Moreover, in view of (\ref{psi-2}), $\Psi(\Omega)$ is equi-continuity. Together with the Arzela-Ascoli theorem, this yields the relative compactness of $\Psi(\Omega)$. Consequently, $\Psi$ is completely continuous. In particular, if $\Omega=\{y\}$ with $y \in C(J,R^n)$, then $\Psi(y)$ is relatively compact. Applying the closeness and convexity of $G$, similarly to step 2 of Lemma \ref{phi-Phi}, we have that $\Psi$ has closed and convex values. Hence the values of $\Psi$ are both compact and convex.

{\bf Step 3.} We show that the upper semicontinuity of $\Psi$.

In view of the complete continuity and compact-valuedness of $\Psi$, it remains to check the closedness of its graph. Let $\{y_n\}$ be a sequence with $y_n \rightarrow y^*$, $\psi_n \in \Psi(y_n)$ and $\psi_n \rightarrow \psi^*$. Then there is $h_n \in \mathbb{S}^1_G(y_n)$ such that
\begin{align*}
\psi_n(t) =& \frac{1}{\Gamma(q)} \int_0^t(t-\tau)^{q-1}g(\tau,y(\tau))h_n(\tau)d\tau- \frac{t}{\Gamma(q)T}\int_0^T(T-\tau)^{q-1}g(\tau,y(\tau))h_n(\tau)d\tau \nonumber \\
&+\frac{t}{T}\int_0^T  c_2(\tau,y_n(\tau))d\tau + \left(1-\frac{t}{T}\right)\int_0^T c_1(\tau,y_n(\tau))d\tau.
\end{align*}
We must show that there exists $h^* \in \mathbb{S}^1_G(y^*)$ such that
\begin{align*}
\psi^*(t) =& \frac{1}{\Gamma(q)} \int_0^t(t-\tau)^{q-1}g(\tau,y(\tau))h^*(\tau)d\tau- \frac{t}{\Gamma(q)T}\int_0^T(T-\tau)^{q-1}g(\tau,y(\tau))h^*(\tau)d\tau \nonumber \\
&+\frac{t}{T}\int_0^T  c_2(\tau,y^*(\tau))d\tau + \left(1-\frac{t}{T}\right)\int_0^T c_1(\tau,y^*(\tau))d\tau.
\end{align*}
Consider the continuous linear operator $\Theta: L^1(J,R^n) \rightarrow C(J,R^n)$ defined by
\begin{align*}
(\Theta h)(t) =& \frac{1}{\Gamma(q)} \int_0^t(t-\tau)^{q-1}g(\tau,y(\tau))h(\tau)d\tau- \frac{t}{\Gamma(q)T}\int_0^T(T-\tau)^{q-1}g(\tau,y(\tau))h(\tau)d\tau \nonumber \\
&+\frac{t}{T}\int_0^T  c_2(\tau,y(\tau))d\tau + \left(1-\frac{t}{T}\right)\int_0^T c_1(\tau,y(\tau))d\tau.
\end{align*}
The definition of $\Theta$ implies that $\psi_n \in \Theta \circ \mathbb{S}^1_G(y_n)$. Similarly to the proof of step 4 in \cite[Theorem 3.1]{WWHWW}, it follows that $\Theta \circ \mathbb{S}^1_G$ has a closed graph. Hence there exists $h^* \in \mathbb{S}^1_G(y^*)$ such that $\Theta(h^*)=\psi^*$. \hfill $\Box$
\end{lemma}

Based on the set-valued Krasnoselskii fixed point theorem, we now state our main result concerning the solvability of BFFDVI (\ref{BFFDVI}).

\begin{theorem}\label{exist-BFFDVI}
Assume that all the conditions of Lemmas \ref{phi-Phi} and \ref{psi-Psi} are satisfied. Then the solution set of BFFDVI (\ref{BFFDVI}) is nonempty.

\textbf{Proof. }\hspace{0.2cm}  According to Remark \ref{Def-Rem}, the existence of mild trajectories of BFFDVI (\ref{BFFDVI}) is equivalent to prove that $\Phi + \Psi$ has a fixed point, where $\Phi$ and  $\Psi$ are defined by (\ref{Phi}) and (\ref{Psi}), respectively. By Lemmas \ref{phi-Phi} and \ref{psi-Psi}, we only need to show that the conclusion (ii) of Lemma \ref{K-fix} is not possible. Let
\begin{equation}\label{delta}
  \delta = \frac{2\frac{M_0+\eta_g \eta_S(1+\eta_Q)}{\Gamma(q+1)}T^q  + (M_1+M_2)T }{1-\rho} +1,
\end{equation}
where $M_0=\underset{t \in J} \sup \left\|\widetilde{F}(t,0)\right\|$. For $\kappa>1$, if there exists $y \in C(J,R^n)$ such that $\kappa y \in \Phi(y) + \Psi(y)$ with $\|y\|_{C}=\delta$. Then there exists $f \in \mathbb{S}^1_{\widetilde{F}}(y)$ and $h \in \mathbb{S}^1_G(y)$ such that
\begin{align*}
\kappa y(t) =& \frac{1}{\Gamma(q)} \int_0^t(t-\tau)^{q-1}\left[f(\tau) + g(\tau,y(\tau))h(\tau)\right]d\tau- \frac{t}{\Gamma(q)T}\int_0^T(T-\tau)^{q-1}\left[f(\tau) + g(\tau,y(\tau))h(\tau)\right]d\tau\\
& +\frac{t}{T}\int_0^T c_2(\tau,y(\tau))d\tau + \left(1-\frac{t}{T}\right)\int_0^T c_1(\tau,y(\tau))d\tau,
\end{align*}
Using (\ref{xFy}) and Remark \ref{G-usc}, one obtains
\begin{eqnarray*}
&&\left\|y(t)\right\|  \nonumber\\
&\leq& \left\|\kappa y(t)\right\|   \nonumber\\
&\leq& \frac{1}{\Gamma(q)} \int_0^t(t-\tau)^{q-1}\left[\|f(\tau)\| + \|g(\tau,y(\tau))\|\|h(\tau)\|\right]d\tau  \\
&&+ \frac{1}{\Gamma(q)}\int_0^T(T-\tau)^{q-1}\left[\|f(\tau)\| + \|g(\tau,y(\tau))\|\|h(\tau)\|\right]d\tau \nonumber\\
&& + \int_0^T \|c_2(\tau,y(\tau))\|d\tau + \int_0^T \|c_1(\tau,y(\tau))\|d\tau \nonumber\\
&\leq&\frac{1}{\Gamma(q)} \int_0^t(t-\tau)^{q-1}\left[L_F \|y(\tau)\| + \left\|\widetilde{F}(\tau,0)\right\| \right]d\tau + \frac{\eta_g \eta_S(1+\eta_Q)}{\Gamma(q)} \int_0^t(t-\tau)^{q-1}d\tau \nonumber\\
&&+\frac{1}{\Gamma(q)} \int_0^T(T-\tau)^{q-1}\left[L_F \|y(\tau)\| + \left\|\widetilde{F}(\tau,0)\right\| \right]d\tau + \frac{\eta_g \eta_S(1+\eta_Q)}{\Gamma(q)} \int_0^T(T-\tau)^{q-1}d\tau\nonumber\\
&& +(M_1+M_2)T \nonumber\\
&\leq&  \frac{2 L_F T^q}{\Gamma(q+1)} \|y\|_C + 2 \frac{M_0+\eta_g \eta_S(1+\eta_Q)}{\Gamma(q+1)} T^q +(M_1+M_2)T,
\end{eqnarray*}
where $M_0=\underset{t \in J} \sup \left\|\widetilde{F}(t,0)\right\|$. This gives 
$$\|y\|_{C} \leq \frac{2 L_F T^q}{\Gamma(q+1)} \|y\|_C + 2 \frac{M_0+\eta_g \eta_S(1+\eta_Q)}{\Gamma(q+1)} T^q +(M_1+M_2)T.$$
Noting that $\|y\|_{C}=\delta$, $\rho= \frac{2 L_F T^q}{\Gamma(q+1)} <1$, we have
$$\delta \leq \frac{2\frac{M_0+\eta_g \eta_S(1+\eta_Q)}{\Gamma(q+1)}T^q  + (M_1+M_2)T }{1-\rho},$$
which contradicts to (\ref{delta}). Consequently, $\Phi + \Psi$ has a fixed point. \hfill $\Box$
\end{theorem}

\section{Stability results for BFFDVI}
\noindent \setcounter{equation}{0}

It is widely recognized that establishing solution stability constitutes a crucial aspect in the well-posedness analysis of DVIs, following the confirmation of solution existence (see, e.g.,\cite{WLZX,WLLH,GLXM,LL2018}). This section concerns the stability analysis of BFFDVI (\ref{BFFDVI}) via the following parametric form, denoted by PBFFDVI:
\begin{equation}\label{PBFFDVI}
\left\{
\begin{array}{l}
\tensor*[^{C}_0]{D}{_t^q} y(t) \in \left[F_{(t,y(t))}\right]_\alpha +g(t,y(t))u(t), \quad  a.e. \; t \in J,\\
u(t) \in \mbox{SOL}\left(\widehat{K}(\varsigma),Q(t,y(t))+\widehat{S}(\cdot, \xi)\right), \quad  a.e. \; t \in J, \\
y(0)=\int_0^T c_1(\tau,y(\tau))d\tau,\; y(T)=\int_0^T c_2(\tau,y(\tau))d\tau,
\end{array}
\right.
\end{equation}
where $\xi$ and $\varsigma$ are two parameters belonging to the metric spaces $(X_1, d_1)$ and $(X_2, d_2)$, respectively. The introduction of parameters $\xi$ and $\varsigma$ perturbs both the mapping $S$ and the set $K$ in BFFDVI (\ref{BFFDVI}), yielding modified mappings $\widehat{S}: R^m \times X_1 \to R^m$ and $\widehat{K}:  X_2 \to 2^{R^m}$. It is worth noting that the stability of BFFDVI (\ref{BFFDVI}) studied here means the semicontinuity properties of the PBFFDVI (\ref{PBFFDVI}) solution map, which helps us predict the trajectory variation of BFFDVI (\ref{BFFDVI}) in response to parametric effects, yielding valuable implications for the design of dynamic systems.

First,  we provide a new existence result for the PBFFDVI  (\ref{PBFFDVI}) based on the set-valued contraction principle. In particular, we present the following hypotheses:
\begin{enumerate}
  \item[(A$_4'$)] there exists $\eta_g, \, L_g>0$ such that $\|g(t,z)\| \leq \eta_g$ for all $(t,z) \in J \times R^n$ and $g(t,\cdot)$ is $L_g$-Lipschitz uniformly for all $t \in J$;
      
  \item[(A$_5'$)] there exists $L_Q>0$ such that the function $Q$ is $L_Q$-Lipschitz, that is,
$$
\|Q(t_1, z_1) - Q(t_2, z_2)\| \leq L_Q (|t_1 - t_2| + \|z_1 - z_2\|)
$$
for all $t_1, t_2 \in J$, and $z_1, z_2 \in R^n$;

  \item[(A$_6'$)] there exists $m_{\widehat{S}}>0$ such that strong monotonicity holds for $\widehat{S}$ with respect to its first argument, that is,
$$
\left<\widehat{S}(u_1, \xi) - \widehat{S}(u_2, \xi), u_1 - u_2 \right> \geq m_{\widehat{S}} \|u_1 - u_2\|^2
$$
for all $u_1, u_2 \in R^m$ and $\xi \in X_2$.
\end{enumerate}

\begin{lemma}\cite[Lemma 3.2]{WLZX}\label{sol_PVI}
Assume that hypotheses (A$_5'$) and (A$_6'$) are satisfied and each value of $\widehat{K}$ is a nonempty convex and closed set. Then, for any given $t\in J$, $y\in C(J,R^n)$, $(\xi,\varsigma) \in X_1 \times X_2$, there exists a unique solution $u(t) \in \widehat{K}(\varsigma)$ to $\mbox{VI}(\widehat{K}(\varsigma), Q(t,y(t))+\widehat{S}(\cdot,\xi))$. Moreover, for fixed $(\xi,\varsigma) \in X_1 \times X_2$ and $y\in C(J,R^n)$, the solution function $u$, regarded as a mapping from the interval $J$ to the set $\widehat{K}(\varsigma)$, is continuous. Furthermore, for $i=1,2$, let $u_i(t)$ be the unique solution to $\mbox{VI}(\widehat{K}(\varsigma), Q(t,y(t))+\widehat{S}(\cdot,\xi))$ associated with $y_i \in C(J,R^n)$ $(i=1,2)$. Then the following inequality holds for all $t \in J$:
\begin{equation}\label{LGQ}
\|u_1(t)-u_2(t)\| \leq \frac{L_Q}{m_{\widehat{S}}}\|y_1(t)-y_2(t)\|.
\end{equation}
\end{lemma}

\begin{remark}
If all assumptions of Lemma \ref{sol_PVI} hold, then $u(t) \in \mbox{SOL}\left(\widehat{K}(\varsigma),Q(t,y(t))+\widehat{S}(\cdot, \xi)\right), \;  a.e. \; t \in J$ in (\ref{PBFFDVI}) can be improved to $u(t) = \mbox{SOL}\left(\widehat{K}(\varsigma),Q(t,y(t))+\widehat{S}(\cdot, \xi)\right)$ for all $t \in J$. Using this fact, we define a new mapping $\widehat{\Upsilon}: C(J,R^n) \rightarrow C(J,R^m)$ as follows
\begin{eqnarray}\label{Upsi}
 \widehat{\Upsilon}(y)=u_y,
 \end{eqnarray}
where $u_y\in C(J,R^m)$ is the unique solution of $\text{VI}(\widehat{K}(\varsigma),Q(t,y(t))+\widehat{S}(\cdot, \xi))$ for each $y \in C(J,R^n)$. 
\end{remark}

\begin{lemma}\label{nonempty-P}
Let (A$_1$)-(A$_3$) and (A$_4'$)-(A$_6'$) hold. Assume that $\widehat{K}$ has nonempty compact and convex values. If there exists $M_{i}>0$ such that $\|c_i(t,z)\|\leq M_{i}$ and there exists $L_{c_i} >0$ such that the function $c_i(t,\cdot)$ is $L_{c_i}$-Lipschitz for all $t\in J$, $z\in R^n$, $i=1,2$, then the function $\widetilde{\varphi}: J \to R^n$ defined by
\begin{align*}
\widetilde{\varphi}(t)
= &\frac{1}{\Gamma(q)}\int_0^t (t-\tau)^{q-1}\chi(\tau)d\tau  -\frac{t}{\Gamma(q)T} \int_0^T (T-\tau)^{q-1} \chi(\tau)d\tau \\
&+ \frac{t}{T} \int_0^T c_2(\tau,y(\tau))d\tau +\left(1-\frac{t}{T} \right)\int_0^T c_1(\tau,y(\tau))d\tau
\end{align*}
is continuous, where $\chi(\tau)=f(\tau)+g(\tau,y(\tau))(\widehat{\Upsilon} y)(\tau)$, $f \in \mathbb{S}^1_{\widetilde{F}}(y)$, $y \in C(J,R^n)$ and $\widehat{\Upsilon}$ is a mapping defined by (\ref{Upsi}).
Moreover, let
\begin{equation}\label{w_lambda}
 \widetilde{\lambda}=2\frac{L_F + L_g \eta_{\widehat{K}}+ \eta_g \frac{L_Q}{m_{\widehat{S}}}}{\Gamma(q+1)}T^q + L_{c_2}T+L_{c_1}T <1,
\end{equation}
where $\eta_{\widehat{K}}$ is defined in (\ref{Upsi_u}).Then the set-valued mapping $\widetilde{\Phi}: C(J,R^n) \to 2^{C(J,R^n)}$ defined by
	\begin{align*}
		\widetilde{\Phi}(y)=& \Bigg\{\widetilde{\varphi} \in C(J, R^n): \widetilde{\varphi}(t)=\frac{1}{\Gamma(q)} \int_0^t(t-\tau)^{q-1}\chi(\tau) d \tau -\frac{t}{\Gamma(q) T} \int_0^T(T-\tau)^{q-1}\chi(\tau) d \tau \\
&\qquad\qquad\qquad\qquad \left.+\frac{t}{T} \int_0^T c_2(\tau, y(\tau)) d \tau+\left(1-\frac{t}{T}\right) \int_0^T c_1(\tau, y(\tau)) d \tau,  f \in \mathbb{S}^1_{\widetilde{F}}(y)\right\}
	\end{align*}
has compact values, and enjoys the contractive property. Furthermore, the nonemptiness of the solution set for PBFFDVI (\ref{PBFFDVI}) is guaranteed.

\textbf{Proof. }\hspace{0.2cm} For any given $(\xi, \varsigma) \in X_1 \times X_2$, the proof is completed in four steps as follows.

{\bf Step 1.} We claim that $\widetilde{\varphi}$ is continuous.

In light of Lemma \ref{FG-nonempty}, we have that $\widetilde{\varphi}$ is well defined. It follows from the compactness of $\widehat{K}(\varsigma)$ and (\ref{Upsi}) that there is a constant $\eta_{\widehat{K}}> 0$ such that
\begin{equation}\label{Upsi_u}
\|\widehat{\Upsilon}(y)\|=\left\|u_y\right\| \leqslant \eta_{\widehat{K}}.
\end{equation}
Given $y \in C(J, R^n )$, $f \in \mathbb{S}^1_{\widetilde{F}}(y)$. In view of the hypothesis (A$_4'$), Remark \ref{Fpr} and (\ref{Upsi_u}), for a.e. $\tau \in J$, one has
\begin{equation}\label{chi_tau}
	\|\chi(\tau)\|=\| f(\tau)+g(\tau, y(\tau))(\widehat{\Upsilon} y)(\tau)\|\leq\| f(\tau)\| + \|g(\tau,y(\tau))\| \|(\widehat{\Upsilon} y)(\tau)\| \leq p(\tau)+\eta_g \eta_{\widehat{K}},
\end{equation}
where $p \in L^\infty(J,R^+)$. Let $0 \leq t_1<t_2 \leq T$,  it follows that
\begin{align}\label{w_phi-1}
\widetilde{\varphi}\left(t_2\right)-\widetilde{\varphi}\left(t_1\right)
 =& \frac{1}{\Gamma(q)} \int_{t_1}^{t_2}(t_2-\tau)^{q-1}\chi(\tau) d\tau +\frac{1}{\Gamma(q)} \int_0^{t_1}\left[(t_2-\tau)^{q-1}-(t_1-\tau)^{q-1}\right]\chi(\tau) d\tau\nonumber \\
& - \frac{t_2-t_1}{\Gamma(q)T}\int_0^T(T-\tau)^{q-1}\chi(\tau)d\tau+ \frac{t_2-t_1}{T}\int_0^T  \left[c_2(\tau,y(\tau))-c_1(\tau,y(\tau)) \right]d\tau.
\end{align}
By (\ref{chi_tau}) and (\ref{w_phi-1}), one gets
\begin{eqnarray}\label{w_psi-2}
 && \|\widetilde{\varphi}\left(t_2\right)-\widetilde{\varphi}\left(t_1\right)\| \nonumber \\
 &\leq&  \frac{1}{\Gamma(q)} \int_{t_1}^{t_2}(t_2-\tau)^{q-1}\|\chi(\tau)\|d\tau +\frac{1}{\Gamma(q)} \int_0^{t_1}\left|(t_2-\tau)^{q-1} -(t_1-\tau)^{q-1}\right|\|\chi(\tau)\|d\tau \nonumber \\
  && + \frac{t_2-t_1}{\Gamma(q)T}\int_0^T(T-\tau)^{q-1} \|\chi(\tau))\|d\tau + \frac{t_2-t_1}{T}\int_0^T  \left\|c_2(\tau,y(\tau))-c_1(\tau,y(\tau)) \right\|d\tau\nonumber \\
   &\leq&   \frac{\|p\|_{L^\infty}+\eta_g \eta_{\widehat{K}}}{\Gamma(q)} \int_{t_1}^{t_2}(t_2-\tau)^{q-1} d\tau +\frac{\|p\|_{L^\infty}+\eta_g \eta_{\widehat{K}}}{\Gamma(q)} \int_0^{t_1}\left[(t_2-\tau)^{q-1}-(t_1-\tau)^{q-1}\right] d\tau\nonumber \\
  &&+ \frac{(t_2-t_1)(\|p\|_{L^\infty}+\eta_g \eta_{\widehat{K}})}{\Gamma(q)T} \int_0^T(T-\tau)^{q-1}d\tau + (M_1+M_2)(t_2-t_1) \nonumber \\
 &=&\left(\frac{(\|p\|_{L^\infty}+\eta_g \eta_{\widehat{K}})T^{q-1}}{\Gamma(q+1)}+M_1+M_2\right) (t_2-t_1) + \frac{\|p\|_{L^\infty}+\eta_g \eta_{\widehat{K}}}{\Gamma(q+1)} \left(t^q_2-t^q_1\right).
\end{eqnarray}
Hence  $\widetilde{\varphi}$ is continuous.

{\bf Step 2.} We show that the compactness of the values of $\widetilde{\Phi}$ holds.

For any $y \in C(J, R^n )$, $\widetilde{\varphi} \in \widetilde{\Phi}(y)$, in light of (\ref{w_psi-2}), one gets the equi-continuity of $\widetilde{\Phi}(y)$. We now proceed to establish the uniform boundedness of $\widetilde{\Phi}(y)$. For each $\widetilde{\varphi} \in \widetilde{\Phi}(y)$, one has
\begin{align*}
	\widetilde{\varphi}(t) =&\frac{1}{\Gamma(q)} \int_0^t(t-\tau)^{\varepsilon-1} \chi(\tau) d \tau-\frac{t}{\Gamma(q) T} \int_0^T(T-\tau)^{q-1} \chi(\tau) d \tau \\
	& +\frac{t}{T} \int_0^T c_2(\tau, y(\tau)) d \tau+\left(1-\frac{t}{T}\right) \int_0^T c_1(\tau, y(\tau)) d \tau .
\end{align*}
It follows from (\ref{chi_tau}) that
\begin{eqnarray*}
 \|\widetilde{\varphi}(t) \|
  &\leq& \frac{1}{\Gamma(q)} \int_0^t(t-\tau)^{q-1}\|\chi(\tau))\|d\tau+ \frac{1}{\Gamma(q)}\int_0^T(T-\tau)^{q-1}\|\chi(\tau))\|d\tau \nonumber \\
&&+ \int_0^T  \|c_2(\tau,y(\tau))\|d\tau + \int_0^T \|c_1(\tau,y(\tau))\|d\tau \nonumber \\
   &\leq& \frac{\|p\|_{L^\infty}+\eta_g \eta_{\widehat{K}}}{\Gamma(q)} \int_0^t(t-\tau)^{q-1} d\tau + \frac{\|p\|_{L^\infty}+\eta_g \eta_{\widehat{K}}}{\Gamma(q)}\int_0^T(T-\tau)^{q-1}d\tau+(M_1+M_2)T \nonumber \\
 &\leq& 2\frac{\|p\|_{L^\infty}+\eta_g \eta_{\widehat{K}}}{\Gamma(q+1)} T^q+(M_1+M_2)T,
\end{eqnarray*}
for all $t\in J$, which implies that
\begin{equation}\label{phiM1M2}
	\left\|\varphi\right\|_C \le 2\frac{\|p\|_{L^\infty}+\eta_g \eta_{\widehat{K}}}{\Gamma(q+1)} T^q+(M_1+M_2)T.
\end{equation}
This gives the uniform boundedness of $\widetilde{\Phi}(y)$. By the Arzela-Ascoli theorem, the relative compactness follows. The closedness of the values of $\widetilde{\Phi}$ is obtained as in step 2 of Lemma \ref{phi-Phi}. Hence, together with the relative compactness, this gives the compactness of the values of $\widetilde{\Phi}$.

{\bf Step 3.}  We show the contractivity of $\widetilde{\Phi}$.

For any $y_1,y_2 \in C(J,R^n)$ and $\widetilde{\varphi}_1 \in \widetilde{\Phi}(y_1)$, there is $\widetilde{f}_1 \in \mathbb{S}^1_{\widetilde{F}}(y_1)$ such that
\begin{align}\label{w_varphi1}
\widetilde{\varphi}_1(t) =& \frac{1}{\Gamma(q)} \int_0^t(t-\tau)^{q-1}\left[\widetilde{f}_1(t)+g(\tau,y_1(\tau))(\widehat{\Upsilon} y_1)(\tau)\right] d\tau  \nonumber \\
&- \frac{t}{\Gamma(q)T}\int_0^T(T-\tau)^{q-1}\left[\widetilde{f}_1(t) + g(\tau,y_1(\tau))(\widehat{\Upsilon} y_1)(\tau)\right]d\tau \nonumber\\
&+\frac{t}{T}\int_{0}^{T}c_2(\tau, y_1(\tau))d\tau+\left(1-\frac{t}{T}\right)\int_{0}^{T}c_1(\tau,y_1(\tau))d\tau.
\end{align}
By the same argument as in step 3 of Lemma \ref{phi-Phi}, we have that there exists an integrable selection $\widetilde{f}_2 \in \mathbb{S}^1_{\widetilde{F}}(y_2)$ such that (\ref{f12F}) holds and $\widetilde{\varphi}_2 \in \widetilde{\Phi}(y_2)$, where
\begin{align}\label{w_varphi2}
\widetilde{\varphi}_2(t) =& \frac{1}{\Gamma(q)} \int_0^t(t-\tau)^{q-1}\left[\widetilde{f}_2(t)+g(\tau,y_2(\tau))(\widehat{\Upsilon} y_2)(\tau)\right] d\tau  \nonumber \\
&- \frac{t}{\Gamma(q)T}\int_0^T(T-\tau)^{q-1}\left[\widetilde{f}_2(t) + g(\tau,y_2(\tau))(\widehat{\Upsilon} y_2)(\tau)\right]d\tau \nonumber\\
&+\frac{t}{T}\int_{0}^{T}c_2(\tau, y_2(\tau))d\tau+\left(1-\frac{t}{T}\right)\int_{0}^{T}c_1(\tau,y_2(\tau))d\tau.
\end{align}
It follows from (\ref{w_varphi1}), (\ref{w_varphi2}) and the Lipschitz condition of $c_i(t,\cdot)$, $i=1,2$, that
\begin{eqnarray}\label{w_varphi} 
&&\left\|\widetilde{\varphi}_1(t)-\widetilde{\varphi}_2(t)\right\| \notag\\
&\le&\frac{1}{\Gamma(q)}\int_{0}^{t}(t-\tau)^{q-1}\left\|\widetilde{f}_1(\tau) +g(\tau,y_1(\tau))(\widehat{\Upsilon} y_1)(\tau)-\widetilde{f}_2(\tau) -g(\tau, y_2(\tau))(\widehat{\Upsilon} y_2)(\tau)\right\|d\tau \notag \\
&&\quad +\frac{1}{\Gamma(q)}\int_{0}^{T}(T-\tau)^{q-1}\left\|\widetilde{f}_1(\tau) +g(\tau, y_1(\tau))(\widehat{\Upsilon} y_1)(\tau)-\widetilde{f}_2(\tau)-g(\tau,y_2(\tau))(\widehat{\Upsilon} y_2)(\tau)\right\|d\tau \notag \\
&&\quad +\int_{0}^{T}\left\|c_2(\tau, y_1(\tau))-c_2(\tau, y_2(\tau))\right\|d\tau + \int_{0}^{T}\left\|c_1(\tau,y_1(\tau))-c_1(\tau,y_2(\tau))\right\|d\tau \notag \\
&\le& \frac{1}{\Gamma(q)}\int_{0}^{t}(t-\tau)^{q-1} \left[\left\|\widetilde{f}_1(\tau)-\widetilde{f}_2(\tau)\right\| + \left\|g(\tau,y_1(\tau))(\widehat{\Upsilon} y_1)(\tau) -g(\tau,y_2(\tau))(\widehat{\Upsilon} y_2)(\tau)\right\|\right]d\tau  \nonumber\\
&&+\frac{1}{\Gamma(q)}\int_{0}^{T}(T-\tau)^{q-1} \left[\left\|\widetilde{f}_1(\tau) -\widetilde{f}_2(\tau)\right\| + \left\|g(\tau,y_1(\tau))(\widehat{\Upsilon} y_1)(\tau) -g(\tau,y_2(\tau))(\widehat{\Upsilon} y_2)(\tau)\right\|\right]d\tau \notag \\
&&+ L_{c_2} T \left\|y_1-y_2\right\|_C + L_{c_1} T\left\|y_1-y_2\right\|_C.
\end{eqnarray}
Using (\ref{LGQ}), (\ref{Upsi}), (\ref{Upsi_u}) and the hypothesis (A$_4'$), we have
\begin{eqnarray}\label{g_tau_Upsi}
&& \left\|g(\tau,y_1(\tau))(\widehat{\Upsilon} y_1)(\tau) -g(\tau,y_2(\tau))(\widehat{\Upsilon} y_2)(\tau)\right\|  \nonumber\\
&=&  \left\|g(\tau,y_1(\tau))(\widehat{\Upsilon} y_1)(\tau) - g(\tau,y_2(\tau))(\widehat{\Upsilon} y_1)(\tau) + g(\tau,y_2(\tau))(\widehat{\Upsilon} y_1)(\tau) -g(\tau,y_2(\tau))(\widehat{\Upsilon} y_2)(\tau)\right\|  \nonumber\\
&\leq&  \left\|g(\tau,y_1(\tau))(\widehat{\Upsilon} y_1)(\tau) - g(\tau,y_2(\tau))(\widehat{\Upsilon} y_1)(\tau)\right\| + \left\|g(\tau,y_2(\tau))(\widehat{\Upsilon} y_1)(\tau) -g(\tau,y_2(\tau))(\widehat{\Upsilon} y_2)(\tau)\right\| \nonumber\\
&\leq&  \left\|g(\tau,y_1(\tau)) - g(\tau,y_2(\tau))\right\|\left\|(\widehat{\Upsilon} y_1)(\tau)\right\| + \left\|g(\tau,y_2(\tau))\right\|\left\|(\widehat{\Upsilon} y_1)(\tau) -(\widehat{\Upsilon} y_2)(\tau)\right\| \nonumber\\
&\leq&  \left(L_g \eta_{\widehat{K}}+ \eta_g \frac{L_Q}{m_{\widehat{S}}}\right)\|y_1(\tau)-y_2(\tau)\|.
\end{eqnarray}
Combining (\ref{f12F}), (\ref{w_varphi}) and (\ref{g_tau_Upsi}),  we obtain
\begin{eqnarray*} 
&&\left\|\widetilde{\varphi}_1(t)-\widetilde{\varphi}_2(t)\right\| \notag\\
&\le& \frac{1}{\Gamma(q)}\int_{0}^{t}(t-\tau)^{q-1} \left[L_F \left\|y_1(\tau)-y_2(\tau)\right\|  + \left(L_g \eta_{\widehat{K}}+ \eta_g \frac{L_Q}{m_{\widehat{S}}}\right)\|y_1(\tau)-y_2(\tau)\|\right]d\tau  \nonumber\\
&&+\frac{1}{\Gamma(q)}\int_{0}^{T}(T-\tau)^{q-1}\left[L_F \left\|y_1(\tau)-y_2(\tau)\right\|  + \left(L_g \eta_{\widehat{K}}+ \eta_g \frac{L_Q}{m_{\widehat{S}}}\right)\|y_1(\tau)-y_2(\tau)\|\right]d\tau \notag \\
&&+ L_{c_2} T \left\|y_1-y_2\right\|_C + L_{c_1} T\left\|y_1-y_2\right\|_C \notag \\
&\le& \frac{L_F + L_g \eta_{\widehat{K}}+ \eta_g \frac{L_Q}{m_{\widehat{S}}}}{\Gamma(q)}\|y_1-y_2\|_C \left(\int_{0}^{t}(t-\tau)^{q-1} d\tau + \int_{0}^{T}(T-\tau)^{q-1} d\tau\right) \nonumber\\
&&+ L_{c_2} T \left\|y_1-y_2\right\|_C + L_{c_1} T\left\|y_1-y_2\right\|_C \nonumber\\
&\le& \left(2\frac{L_F + L_g \eta_{\widehat{K}}+ \eta_g \frac{L_Q}{m_{\widehat{S}}}}{\Gamma(q+1)}T^q + L_{c_2}T+L_{c_1}T \right) \|y_1-y_2\|_C = \widetilde{\lambda} \|y_1-y_2\|_C,
\end{eqnarray*}
where $\widetilde{\lambda}$ is defined by (\ref{w_lambda}), and so
\begin{equation*}
	\left\|\widetilde{\varphi}_1-\widetilde{\varphi}_2\right\|_C \le \widetilde{\lambda} \left\|y_1-y_2\right\|_C.
\end{equation*}
It follows that \begin{equation*}
d(\widetilde{\varphi}_1,\widetilde{\Phi}(y_2)) = \underset{\widetilde{\varphi}_2 \in \widetilde{\Phi}(y_2)}\inf\| \widetilde{\varphi}_1-\widetilde{\varphi}_2\|_{C}
\leq \widetilde{\lambda} \left\|y_1-y_2\right\|_{C}.
\end{equation*}
Since $\widetilde{\varphi}_1\in \widetilde{\Phi}(y_1)$ is arbitrary, we have
\begin{equation*}
\underset{\widetilde{\varphi}_1\in \widetilde{\Phi}(y_1)}\sup
d(\widetilde{\varphi}_1,\widetilde{\Phi}(y_2)) \leq \widetilde{\lambda} \left\|y_1-y_2\right\|_{C}.
\end{equation*}
Similarly, one has
\begin{equation*}
\underset{\widetilde{\varphi}_2 \in \widetilde{\Phi}(y_2)}\sup d(\widetilde{\varphi}_2,\widetilde{\Phi}(y_1)) \leq \widetilde{\lambda} \left\|y_1-y_2\right\|_{C}.
\end{equation*}
Applying the definition of the Hausdorff distance, one has
\begin{equation*}
H\left(\widetilde{\Phi}(y_1),\widetilde{\Phi}(y_2)\right) \leq \widetilde{\lambda} \left\|y_1-y_2\right\|_{C},
\end{equation*}
which implies that $\widetilde{\Phi}$ is contractive since $\widetilde{\lambda}<1$.

{\bf Step 4.}  We show that the solution set of PBFFDVI (\ref{PBFFDVI}) is nonempty.

According to Remark \ref{Def-Rem}, Mild trajectories of PBFFDVI (\ref{PBFFDVI}) exist if and only if $\widetilde{\Phi}$ has a fixed point, which is guaranteed by Lemma \ref{Fixpoint}. \hfill $\Box$
\end{lemma}

Next, we show that the mapping $(\varsigma, \xi) \mapsto \text{SOL}(\text{PBFFDVI}(\varsigma,\xi))$ is closed, where $\text{SOL}(\text{PBFFDVI}(\varsigma,\xi))$ stands for the solution set of PBFFDVI (\ref{PBFFDVI}).

\begin{lemma}\label{SOLcl}
Let all assumptions of Lemma \ref{nonempty-P} be satisfied. If the function $\widehat{S}$ is continuous and the set-valued mapping $\widehat{K}$ is continuous, then $\text{SOL}(\text{PBFFDVI}(\cdot,\cdot))$ is closed.

\textbf{Proof.}\hspace{0.2cm} Consider any $(\xi, \varsigma) \in X_1 \times X_2$. By Lemma \ref{nonempty-P}, $\text{SOL}(\text{PBFFDVI}(\xi, \varsigma)) \neq \emptyset$. To establish the desired closedness of the set-valued mapping $ \text{SOL}(\text{PBFFDVI}(\cdot,\cdot))$, it suffices to prove that the graph  $Gr(\text{SOL}(\text{PBFFDVI}(\cdot,\cdot)))=\{((\xi, \varsigma),(y,u)) \in (X_1 \times X_2) \times (R^n\times R^m):  (y,u) \in \text{SOL}(\text{PBFFDVI}(\xi, \varsigma))\}$ is a closed set. It is worth noting that, according to Lemmas \ref{sol_PVI} and \ref{nonempty-P}, $y$ and $u$ are two continuous functions. Let $\{(\xi_n, \varsigma_n)\} \subset X_1 \times X_2$ and $\{(y_n, u_n)\} \subset C(J, R^n) \times C(J,R^m)$ be given such that 
$$(y_n, u_n) \in \mbox{SOL}(\mbox{PBFFDVI}(\xi_n, \varsigma_n)).$$
And assume that $(\xi_n, \varsigma_n) \rightarrow (\xi, \varsigma)$ in $X_1 \times X_2$, $y_n \rightarrow y$ in $C(J, R^n)$ and $u_n \rightarrow u$ in $C(J,R^m)$. Consequently, for every $n\geq1$, there exists
$f_n \in \mathbb{S}^1_{\widetilde{F}}(y_n)$ such that
\begin{equation}\label{xyn}
\left\{
\begin{array}{l}
y_n(t)=w_n(t)+\frac{1}{\Gamma(q)}\int_0^t (t-\tau)^{q-1}g(\tau,y_n(\tau))u_n(\tau)d\tau  -\frac{t}{\Gamma(q)T} \int_0^T (T-\tau)^{q-1} g(\tau,y_n(\tau))u_n(\tau)d\tau\\
\qquad\qquad + \frac{t}{T} \int_0^T c_2(\tau,y_n(\tau))d\tau +\left(1-\frac{t}{T} \right)\int_0^T c_1(\tau,y_n(\tau))d\tau, \; t \in J, \\
u_n(t) \in \widehat{K}(\varsigma_n), \\
\left<Q(t,y_n(t))+\widehat{S}(u_n(t),\xi_n),v-u_n(t)\right> \geq 0, \; \forall \;  v \in \widehat{K}(\varsigma_n), \; t \in J, 
\end{array}
\right.
\end{equation}
where
\begin{equation}\label{yn}
w_n(t)=\frac{1}{\Gamma(q)} \int_0^t(t-\tau)^{q-1}f_n(\tau) d\tau- \frac{t}{\Gamma(q)T}\int_0^T(T-\tau)^{q-1}f_n(\tau)d\tau.
\end{equation}

By the same argument as Lemma \ref{phi-Phi} in \cite{WLZX}, it follows from (\ref{xyn}) that
\begin{eqnarray}
&&u(t) \in \widehat{K}(\varsigma), \nonumber\\
&&\left<Q(t,y(t))+\widehat{S}(u(t),\xi),v-u(t)\right> \geq 0, \; \forall \;  v \in \widehat{K}(\varsigma), \; t \in J. \nonumber
\end{eqnarray}

We recall that $y_n \rightarrow y$ in $C(J, R^n)$ and $u_n \rightarrow u$ in $C(J,R^m)$, which implies $y_n(t) \rightarrow y(t)$ and $u_n(t) \rightarrow u(t)$ for all $t \in J$. Moreover, one can find a constant $M_3>0$ satisfying $\|u_n(t)\| \leq M_3$ for all $n\geq1$ and $t \in J$. Since $g(t,\cdot)$ and $c_i(t,\cdot) \, (i=1,2)$ are Lipschitz for every $t \in J$, the mapping $g(t,\cdot)$ and $c_i(t,\cdot) \, (i=1,2)$ are continuous. Combining this with $y_n(\tau) \rightarrow y(\tau)$ and $u_n(\tau) \rightarrow u(\tau)$, we obtain, for all $\tau \in (0,t]$,
$$(t-\tau)^{q-1}g(\tau,y_n(\tau))u_n(\tau) \rightarrow (t-\tau)^{q-1}g(\tau,y(\tau))u(\tau)$$
and
$$c_i(\tau,y_n(\tau)) \rightarrow c_i(\tau,y(\tau)), i=1,2.$$
In addition, since $g$ and $c_i \, (i=1,2)$ are bounded functions, for all $n$,
$$\|(t-\tau)^{q-1}g(\tau,y_n(\tau))u_n(\tau)\| \leq |t-\tau|^{q-1} \eta_g M_3$$
and
$$\|c_i(\tau,y_n(\tau))\| \leq M_i,i=1,2.$$
Applying the Lebesgue dominated convergence theorem, for each $t \in J$, it follows that
\begin{equation}\label{Bu}
\frac{1}{\Gamma(q)}\int_{0}^{t}(t-\tau)^{q-1}g(\tau,y_n(\tau))u_n(\tau)d\tau \rightarrow \frac{1}{\Gamma(q)}\int_{0}^{t}(t-\tau)^{q-1}g(\tau,y(\tau))u(\tau)d\tau,
\end{equation}
\begin{equation}\label{Cu1}
\left(1-\frac{t}{T} \right)\int_0^T c_1(\tau,y_n(\tau))d\tau \rightarrow  \left(1-\frac{t}{T} \right)\int_0^T c_1(\tau,y(\tau))d\tau
\end{equation}
and
\begin{equation}\label{Cu2}
\frac{t}{T} \int_0^T c_2(\tau,y_n(\tau))d\tau  \rightarrow  \frac{t}{T} \int_0^T c_2(\tau,y(\tau))d\tau.
\end{equation}
As in the proof of closedness in step 2 of Lemma \ref{phi-Phi}, Remark \ref{Fpr} gives a subsequence of $\{f_n\}$ (still denoted by $\{f_n\}$) such that $x_n=\sum\limits_{j=n}^{j_0(n)}\rho_{n,j}f_j$ $(n \geq 1)$ converges to $\widetilde{f}$ in $L^1(J,R^n)$, where $j_0(n)$ and $\rho_{n,j}$ are defined as (\ref{x-rho-f}). If necessary, pass to a further subsequence (again denoted by $\{x_n\}$) to get $x_n(\tau) \rightarrow \widetilde{f}(\tau)$ for $a.e.$ $\tau \in J$. By the same arguments in (4.5) and (4.6) of \cite{WLZX}, it follows that, for $a.e.$ $\tau \in J$,
\begin{eqnarray*} 
\widetilde{f}(\tau) & \in & \mathop{\bigcap}\limits_{n \geq1} \overline{\{x_j(\tau):j\geq n\}}  \subseteq \mathop{\bigcap}\limits_{n \geq1} \overline{\mbox{conv}}\{f_j(\tau): j\geq n\} \notag\\
&\subseteq &  \mathop{\bigcap}\limits_{n \geq1} \overline{\mbox{conv}} \left( \underset{j \geq n}\bigcup \widetilde{F}(\tau,y_j(\tau))\right)  \nonumber\\
&=&  \overline{\mbox{conv}}\left(\limsup \limits_{j \rightarrow \infty} \widetilde{F}(\tau,y_j(\tau))\right) \nonumber\\
&\subseteq & \overline{\mbox{conv}}\left(\widetilde{F}(\tau,y(\tau))\right).
\end{eqnarray*}
Given the compact convex values of $\widetilde{F}$, one gets that  $\widetilde{f}(\tau) \in \widetilde{F}(\tau,y(\tau))$ for $a.e.$ $\tau \in J$. Thus $\widetilde{f} \in \mathbb{S}^1_{\widetilde{F}}(y)$. Let $\widetilde{y}_n=\sum\limits_{j=n}^{j_0(n)}\rho_{n,j}y_j$. It then follows from (\ref{xyn}) and (\ref{yn}) that
\begin{align}\label{xyn1}
  \widetilde{y}_n(t) =& \widetilde{w}_n(t) + \sum\limits_{j=n}^{j_0(n)}\frac{\rho_{n,j}}{\Gamma(q)} \int_0^t (t-\tau)^{q-1}g(\tau,y_j(\tau))u_j(\tau)d\tau  - \frac{t}{T}\sum\limits_{j=n}^{j_0(n)}\frac{ \rho_{n,j} }{\Gamma(q)} \int_0^T (T-\tau)^{q-1} g(\tau,y_j(\tau))u_j(\tau)d\tau \nonumber\\
& + \sum\limits_{j=n}^{j_0(n)}\frac{t \rho_{n,j} }{T} \int_0^T c_2(\tau,y_j(\tau))d\tau + \sum\limits_{j=n}^{j_0(n)}\rho_{n,j}\left(1-\frac{t}{T} \right)\int_0^T c_1(\tau,y_j(\tau))d\tau, \; t \in J, 
\end{align}
where
\begin{equation}\label{yn1}
\widetilde{w}_n(t)=\frac{1}{\Gamma(q)} \int_0^t(t-\tau)^{q-1}x_n(\tau)d\tau - \frac{t}{\Gamma(q)T}\int_0^T(T-\tau)^{q-1}x_n(\tau)d\tau.
\end{equation}
By virtue of $y_n \rightarrow y$, (\ref{Bu}), (\ref{Cu1}) and (\ref{Cu2}), we directly conclude that $\widetilde{y}_n(t) \rightarrow y(t)$ and
\begin{equation*}
\sum\limits_{j=n}^{j_0(n)}\frac{\rho_{n,j}}{\Gamma(q)} \int_0^t (t-\tau)^{q-1}g(\tau,y_j(\tau))u_j(\tau)d\tau \rightarrow \frac{1}{\Gamma(q)}\int_{0}^{t}(t-\tau)^{q-1}g(\tau,y(\tau))u(\tau)d\tau,
\end{equation*}
\begin{equation*}
\sum\limits_{j=n}^{j_0(n)}\rho_{n,j}\left(1-\frac{t}{T} \right)\int_0^T c_1(\tau,y_j(\tau))d\tau \rightarrow  \left(1-\frac{t}{T} \right)\int_0^T c_1(\tau,y(\tau))d\tau.
\end{equation*}
and
\begin{equation*}
\sum\limits_{j=n}^{j_0(n)}\frac{t \rho_{n,j} }{T} \int_0^T c_2(\tau,y_j(\tau))d\tau  \rightarrow  \frac{t}{T} \int_0^T c_2(\tau,y(\tau))d\tau.
\end{equation*}
for all $t \in J$. In light of Remark \ref{Fpr}, for every $t \in J$, almost every $\tau \in (0,t]$ and each $n\geq1$, one has
$$\left\|(t-\tau)^{q-1}x_n(\tau)\right\| \leq  (t-\tau)^{q-1}  \sum\limits_{j=n}^{j_0(n)}\rho_{n,j}\left\|f_j(\tau) \right\|\leq \left(t-\tau\right)^{q-1} p(\tau).$$
It should be noted that $\left(t-\cdot\right)^{q-1} p(\cdot)  \in L^1((0,T],R^+)$ and $x_n(\tau) \rightarrow \widetilde{f}(\tau)$ for $a.e.$ $\tau \in J$. Let $n \rightarrow \infty$ in (\ref{xyn1}) and (\ref{yn1}). Applying the Lebesgue dominated convergence theorem that, we get, for all $t \in J$,
\begin{align*}
  y(t) =& w(t) + \frac{1}{\Gamma(q)} \int_0^t (t-\tau)^{q-1}g(\tau,y(\tau))u(\tau)d\tau  - \frac{t}{\Gamma(q)T} \int_0^T (T-\tau)^{q-1} g(\tau,y(\tau))u(\tau)d\tau \nonumber\\
& + \frac{t}{T} \int_0^T c_2(\tau,y(\tau))d\tau + \left(1-\frac{t}{T} \right)\int_0^T c_1(\tau,y(\tau))d\tau, \; t \in J, 
\end{align*}
where
\begin{equation*}
w(t)=\frac{1}{\Gamma(q)} \int_0^t(t-\tau)^{q-1}\widetilde{f}(\tau)d\tau - \frac{t}{\Gamma(q)T}\int_0^T(T-\tau)^{q-1}\widetilde{f}(\tau)d\tau, \; \widetilde{f} \in \mathbb{S}^1_{\widetilde{F}}(y).
\end{equation*}
Consequently, $(y,u) \in \mbox{SOL}(\mbox{PBFFDVI}(\xi, \varsigma))$, that is, $\mbox{SOL}(\mbox{PBFFDVI}(\cdot, \cdot))$ is closed.  \hfill $\Box$
\end{lemma}

We now present the u.s.c. of the solution mapping of (\ref{PBFFDVI}).

\begin{theorem}\label{usc-PBFFDVI}
Assume that all assumptions of Lemma \ref{SOLcl} are satisfied. If the set-valued mapping $\widehat{K}: X_2 \rightarrow 2^{R^m}$ is locally compact, then $\mbox{SOL}(\mbox{PBFFDVI}(\cdot, \cdot))$ is u.s.c.

\textbf{Proof.}\hspace{0.2cm} Take $(\xi^*, \varsigma^*) \in X_1 \times X_2$ and consider a neighborhood $\Xi_1 \times \Xi_2$ of $(\xi^*, \varsigma^*)$ in $X_1 \times X_2$ satisfying the closure of $\underset{\varsigma \in \Xi_2}\bigcup \widehat{K}(\varsigma)$ is compact in $R^m$. By Lemma  \ref{nonempty-P}, it follows that $\mbox{SOL}(\mbox{PBFFDVI}(\xi, \varsigma))$ is nonempty for every $(\xi^*, \varsigma^*) \in \Xi_1 \times \Xi_2$. Suppose that $\mbox{SOL}(\mbox{PBFFDVI}(\cdot, \cdot))$ fails to be u.s.c. at $(\xi^*, \varsigma^*)$. Then one can find an open set $O$ containing $\mbox{SOL}(\mbox{PBFFDVI}(\xi^*, \varsigma^*))$, a sequence $\{(\xi_n, \varsigma_n)\} \subset \Xi_1 \times \Xi_2$ converging to $(\xi^*, \varsigma^*)$ and $\{(y_n,u_n)\} \subset C(J, R^n) \times C(J,R^m)$ with
$$(y_n,u_n) \in \mbox{SOL}(\mbox{PBFFDVI}(\xi_n, \varsigma_n))$$
such that $(y_n,u_n) \notin O$ for all $n \geq 1$. As $(y_n,u_n) \in \mbox{SOL}(\mbox{PBFFDVI}(\xi_n, \varsigma_n))$, it holds that
\begin{eqnarray*}
&&u_n(t) \in \widehat{K}(\varsigma_n), \nonumber\\
&&\left<Q(t,y_n(t))+\widehat{S}(u_n(t),\xi_n),v-u_n(t)\right> \geq 0, \; \forall \;  v \in \widehat{K}(\varsigma_n),  \; t \in J
\end{eqnarray*}
and
\begin{align*}
y_n(t)=&\frac{1}{\Gamma(q)}\int_0^t (t-\tau)^{q-1}\left[f_n(\tau)+ g(\tau,y_n(\tau))u_n(\tau)\right]d\tau  \\
&- \frac{t}{\Gamma(q)T} \int_0^T (T-\tau)^{q-1} \left[f_n(\tau) + g(\tau,y_n(\tau))u_n(\tau)\right]u_n(\tau)d\tau  \\
&+ \frac{t}{T} \int_0^T c_2(\tau,y_n(\tau))d\tau +\left(1-\frac{t}{T} \right)\int_0^T c_1(\tau,y_n(\tau))d\tau, \; t \in J,
\end{align*}
where $f_n \in \mathbb{S}^1_{\widetilde{F}}(y_n)$. A key observation is recorded below.

{\bf Claim 1.} The sequence $\{y_n\}$ is relatively compact.

By the compactness of $\overline{\underset{\varsigma \in \Xi_2}\bigcup \widehat{K}(\varsigma)}$, there is a constant $M^*>0$ such that
\begin{equation}\label{unM*}
 \|u_n(t)\| \leq  M^*
\end{equation}
for all $t \in J$. Similarly to (\ref{w_psi-2}) and (\ref{phiM1M2}), we obtain
$$\|y_n(t_2) - y_n(t_1)\| \leq \left(\frac{\left(\|p\|_{L^\infty}+\eta_g M^*\right)T^{q-1}}{\Gamma(q+1)}+M_1+M_2\right) (t_2-t_1) + \frac{\|p\|_{L^\infty}+\eta_g M^*}{\Gamma(q+1)} \left(t^q_2-t^q_1\right)$$
for $0 \leq t_1 < t_2 \leq T$, and
$$\|y_n\|_{C} \leq 2\frac{\|p\|_{L^\infty}+\eta_g M^*}{\Gamma(q+1)} T^q+(M_1+M_2)T.$$ 
Thus $\{y_n\}$ is equi-continuous and uniformly bounded. By the Arzela-Ascoli theorem, $\{y_n\}$ is relatively compact.

{\bf Claim 2.} The sequence $\{u_n\}$ is relatively compact.

By virtue of (\ref{unM*}), $\{u_n\}$ is uniformly bounded in $C(J,R^m)$. In addition, similarly to the proof of (3.4) in \cite{WLZX}, one has
\begin{equation*}
\|u_n(t_2)-u_n(t_1)\| \leq \frac{L_Q}{m_{\widehat{S}}}(|t_2-t_1|+ \left\|y_n(t_2)-y_n(t_1)\right\|).
\end{equation*}
Noting that $\{y_n\}$ is equi-continuous by claim 1, we conclude that $\{u_n\}$ is equi-continuous as well. Consequently, by the Arzela-Ascoli theorem, $\{u_n\}$ is relatively compact.

From claims 1 and 2, we choose subsequences, still denoted by $\{y_n\}$ and $\{u_n\}$, with $y_n \rightarrow y^*$ in $C(J, R^n)$ and $u_n \rightarrow u^*$ in $C(J,R^m)$. Lemma \ref{SOLcl} gives the closedness of $\mbox{SOL}(\mbox{PBFFDVI}(\cdot, \cdot))$, so $(y^*,u^*)$ belong to $\mbox{SOL}(\mbox{PBFFDVI}(\xi^*, \varsigma^*))$, which is contained in $O$. On the other hand, since every $(y_n,u_n)$ lines outside $O$, the limit point $(y^*,u^*)$ must stay outside $O$ as well. This is impossible. Consequently, $\mbox{SOL}(\mbox{PBFFDVI}(\cdot, \cdot))$ is u.s.c. at $(\xi^*, \varsigma^*) \in X_1 \times X_2$.   \hfill $\Box$
\end{theorem}

\begin{remark}
Motivated by the modeling of microbial metabolism with memory effects in uncertain environments, Wu et al. \cite{WWHWW} proposed a fuzzy fractional DVI in 2021. This system combines the core features of FFDIs with fractional DVIs. The authors studied the existence of solutions, constructed numerical algorithms, and provided two numerical examples to verify their theoretical results. In 2024, Wu et al. \cite{WLZX} established a new existence result for the FFDVI. Their work also showed the compactness of the solution set, the continuous dependence of solutions on initial values, and the u.s.c. of the solution mapping with respect to parameters. This paper focuses on the BFFDVI. We establish a theoretical framework for fuzzy fractional boundary value problems constrained by variational inequalities. This framework can be employed to explore boundary value problems arising from physical processes and materials with memory effects or anomalous phenomena under uncertainty. We then investigate the solvability and stability of the BFFDVI. 
\end{remark}

\section{Numerical examples}
\noindent \setcounter{equation}{0}
In this section, we illustrate the main results with two numerical examples. In particular, we use $\|\cdot\|_i$ $(i=1,2)$ denotes the $\ell_i$-norm on $R^n$.
\begin{example}
We consider the following BFFDVI
\begin{equation}\label{ex1}
\left\{
\begin{array}{l}
\leftidx{_0^C}D{_t^{1.6}}y(t) \in \cos(y(t)) \cdot [\omega]_\alpha +
\left( \begin{array}{c}
1.2 \sin t \\
-2.5 \cos y(t)\\
\end{array} \right)^\top u(t), \quad  a.e. \; t \in [0,0.7],\\
\left<\left(\begin{array}{c}
\arctan y(t) + 2\pi \\
-1.4 e^{-t}\\
\end{array} \right) + \left( \begin{array}{cc}
    3 & 0\\
    0 & 3 \\
\end{array}
\right) u(t),v-u(t)\right> \geq 0, \; \forall \; v \in K, \; a.e. \; t \in [0,0.7], \\
y(0)=1.2\int_0^{0.7}\sin(y(\tau))d\tau,\; y(0.7)=0.9\int_0^{0.7} \cos(y(\tau))d\tau,
\end{array}
\right.
\end{equation}
where $K=\left\{u=\left(u_{1},u_{2}\right)^{\top}| 0 \leq u_{1} < \infty , \;  0 \leq u_{2}  < \infty \right\}$, $\alpha \in [0,1]$, $u(t)=\left(u_{1}(t),u_{2}(t)\right)^{\top}$, $y(t) \in R$, $\omega$ is a symmetric triangular fuzzy number (STFN), whose level sets are
$$[\omega]_\alpha =[0.5(\alpha-1),0.5(1-\alpha)].$$

Set
$$F_{(t,y)}=\cos (y) \cdot \omega,
\quad
g(t,y)=\left(
         \begin{array}{c}
           1.2 \sin t \\
           -2.5 \cos y(t)\\
         \end{array}
       \right)^\top,
$$
$$
Q(t,y)=\left(
         \begin{array}{c}
          \arctan y + 2\pi \\
           -1.4 e^{-t}\\
         \end{array}
       \right),
       \quad
S(u)=\left(
         \begin{array}{cc}
          3 & 0\\
          0 & 3 \\
         \end{array}
       \right)\left(
         \begin{array}{c}
          u_1 \\
           u_2 \\
         \end{array}
       \right)=\left(
         \begin{array}{c}
          3u_1 \\
          3u_2 \\
         \end{array}
       \right),$$
$$c_1(t,y)=1.2\sin(y),\quad c_2(t,y)=0.9\cos(y).$$
Obviously, $g,Q,S,c_i \ (i=1,2)$ are continuous functions, $Q$ is monotone and assumption (A$_2$) is satisfied. $[\omega]_0=[-0.5,0.5]$. Moreover, for any $y_1, y_2, y \in R$ and any $t \in [0,0.7]$, we have
\begin{eqnarray*}
&&\mathcal{H}\left(F_{(t,y_1)},F_{(t,y_2)}\right)  \nonumber\\
&=& \underset{\alpha \in[0,1]}\sup H\left(\cos (y_1) \cdot [w]_\alpha,\cos (y_2) \cdot [w]_\alpha\right) \nonumber\\
&=&\underset{\alpha \in[0,1]}\sup H\left(\cos (y_1) \cdot \left[0.5(\alpha-1), 0.5(1-\alpha)\right], \cos (y_2) \cdot \left[0.5(\alpha-1), 0.5(1-\alpha)\right]\right) \nonumber\\
&=&\underset{\alpha \in[0,1]}\sup H\left(0.5(1-\alpha)\left[-\left|\cos (y_1)\right|, \left|\cos (y_1)\right|\right], 0.5(1-\alpha)\left[-\left|\cos (y_2)\right|, \left|\cos (y_2)\right| \right]\right) \nonumber\\
&=&\underset{\alpha \in[0,1]}\sup  0.5(1-\alpha)\big|\left|\cos (y_1)\right|-\left|\cos (y_2)\right|\big| \nonumber\\
&\leq&\underset{\alpha \in[0,1]}\sup  0.5(1-\alpha)\big|\cos (y_1)-\cos (y_2)\big| \nonumber\\
&\leq& 0.5 \big|-\sin(\zeta) \cdot \left(y_1-y_2\right) \big|           \nonumber\\
&\leq& 0.5 \left|y_1-y_2 \right| =  L_F \left|y_1-y_2 \right| ,
\end{eqnarray*}
$$\|F_{(t,y)}\| = \mathcal{H}\left(F_{(t,y)},\widetilde{0}\right) = \underset{x \in \left[F_{(t,y)}\right]_0} \sup |x| = \underset{x \in \cos (y) \cdot[-0.5,0.5]} \sup |x| \leq 0.5 =p(t),$$
$$\|g(t,y)\|_1=1.2|\sin t|+2.5|\cos y| \leq 3.7 = \eta_g,$$
$$\|Q(t,y)\|_1=|\arctan y+2\pi|+\left|-1.4 e^{-t}\right| \leq \frac{5\pi}{2} +1.4 = \eta_Q,$$
$$\lim \limits_{u \in K, \|u\|_2 \rightarrow \infty} \frac{\left<S(u),u-u_0\right>}{\|u\|_2^2}=3>0,$$
where $\zeta$ exists in between $y_1$ and $y_2$, $\widetilde{0}$ is a fuzzy set defined by $\widetilde{0}(x)=1$ if $x=0$ and $\widetilde{0}(x)=0$ if $x\neq0$, $u_0=(0,0)^{\top}$. This shows that assumptions (A$_1$) and (A$_3$)-(A$_6$) hold. Moreover, one has
$$|c_1(t,y)| \leq 1.2= M_1, \quad  |c_2(t,y)| \leq 0.9 = M_2,$$
$$\rho=\frac{2 L_F T^q}{\Gamma(q+1)} = 0.3953 < 1$$
for all $t \in [0,T]$ and all $y\in R$. This verifies all the assumptions of Theorem \ref{exist-BFFDVI}. Hence BFFDVI (\ref{ex1}) admits at least one solution.
\end{example}

\begin{example}
We consider the following PBFFDVI
\begin{equation}\label{ex2}
\left\{
\begin{array}{l}
\leftidx{_0^C}D{_t^{1.8}}y(t) \in \sin y(t) \cdot [\omega]_\alpha +
\left( \begin{array}{c}
-0.47 \sin y(t)\\
1.43 \cos t \\
\end{array} \right)^\top u(t), \quad  a.e. \; t \in [0,0.5],\\
\left<\left(\begin{array}{c}
-0.7 \cos y(t) \\
1.2\arctan t\\
\end{array} \right) + \left(
         \begin{array}{c}
          -0.55 \arctan u_1(t) +4.1 u_1(t) +0.02 \xi_1\\
           0.12 \sin u_2(t) +3.72 u_2(t) -0.01 \sin \xi_2\\
         \end{array}
       \right),v-u(t)\right> \geq 0, \\
        \forall \; v \in \widehat{K}(\varsigma), \; a.e. \; t \in [0,0.5], \\
y(0)=0.4\int_0^{0.5}\cos(y(\tau))d\tau,\; y(0.5)=0.3\int_0^{0.5} \arctan(y(\tau))d\tau,
\end{array}
\right.
\end{equation}
where $\widehat{K}(\varsigma)=\left\{\left(u_{1},u_{2}\right)^{\top}\left|u_{1}^{2}+u_{2}^{2} \leq (0.72+0.18 \cos \varsigma)^2 \right.\right\}$, $u(t)=(u_1(t),u_2(t))^\top \in R^2$, $y(t) \in R$, $\varsigma\in R$, $\xi=(\xi_1,\xi_2)^\top \in R^2$, $\alpha \in [0,1]$, $\omega$ is a STFN with its level set as follows
$$[\omega]_\alpha =[0.6(\alpha-1),0.6(1-\alpha)].$$

Set
$$F_{(t,y)}=\sin y \cdot \omega,
\quad
g(t,y)=\left(
         \begin{array}{c}
           -0.47 \sin y\\
           1.43 \cos t\\
         \end{array}
       \right)^\top,
$$
$$
Q(t,y)=\left(
         \begin{array}{c}
          - 0.7 \cos y \\
           1.2\arctan t \\
         \end{array}
       \right),
$$
$$
\widehat{S}(u,\xi)=\left(
         \begin{array}{c}
          -0.55 \arctan u_1(t) +4.1 u_1(t) +0.02 \xi_1\\
           0.12 \sin u_2(t) +3.72 u_2(t) -0.01 \sin \xi_2\\
         \end{array}
       \right),$$
$$c_1(t,y)=0.4\cos(y),\quad c_2(t,y)=0.3\arctan(y).$$
Obviously, $g, Q,\widehat{S}, c_i \ (i=1,2)$ are continuous, and $\widehat{K}$ has nonempty compact convex values and is locally compact. Moreover, assumption (A$_2$) is satisfied. For any $\varsigma_1, \varsigma_2 \in R$, one has
$$H\left(\widehat{K}\left(\varsigma_1\right),\widehat{K}\left(\varsigma_2\right)\right) = \max\left\{\underset{\varrho_1\in \widehat{K}\left(\varsigma_1\right)}\sup \, d\left(\varrho_1,\widehat{K}\left(\varsigma_2\right)\right), \underset{\varrho_2\in \widehat{K}(\varsigma_2)}\sup \, d\left(\varrho_2,\widehat{K}\left(\varsigma_1\right)\right) \right\} \leq 0.18\left|\varsigma_1-\varsigma_2\right|,$$
which means that $\widehat{K}$ is Lipschitz. Thus $\widehat{K}$ is continuous.

In addition, $[\omega]_0=[-0.6,0.6]$, and for any $t_1, t_2, t \in [0,0.5]$, $y_1, y_2, y \in R$, one has
$$\mathcal{H}\left(F_{(t,y_1)},F_{(t,y_2)}\right)\leq 0.6 \left|y_1-y_2\right| = L_F \left|y_1-y_2\right|,$$
$$\|F_{(t,y)}\|  \leq 0.6 =p(t), \;\|g(t,y)\|_2 \leq \|g(t,y)\|_1 \leq 1.9 = \eta_g,$$
$$\|g(t,y_1)-g(t,y_2)\|_2 \leq \|g(t,y_1)-g(t,y_2)\|_1 \leq 0.47 \left|y_1-y_2 \right| = L_g \left|y_1-y_2 \right|,$$
$$|c_1(t,y)| \leq 0.4= M_1, \; |c_2(t,y)| \leq 0.3 = M_2,$$
$$|c_1(t,y_1)-c_1(t,y_2)|\leq 0.4 \left|y_1-y_2 \right| = L_{c_1} \left|y_1-y_2 \right|,$$
$$ |c_2(t,y_1)-c_2(t,y_2)|\leq 0.3 \left|y_1-y_2 \right| = L_{c_2} \left|y_1-y_2 \right|,$$
$$\left\|Q(t_1,y_1)-Q(t_2,y_2) \right\|_2  \leq \left\|Q(t_1,y_1)-Q(t_2,y_2) \right\|_1  \leq 1.2 (|t_1-t_2|+\left|y_1-y_2 \right|) = L_Q (|t_1-t_2|+\left|y_1-y_2 \right|).$$
This shows that assumptions (A$_1$), (A$_3$), (A$'_4$) and (A$'_5$) hold.

Furthermore, for any $\xi \in R^2$, $u^i=(u^i_1,u^i_2)^\top \in R^2$ $(i=1,2)$, we have
\begin{eqnarray*}
&& \left<\widehat{S}\left(u^1,\xi\right) - \widehat{S}\left(u^2,\xi\right),u^1-u^2\right>  \nonumber\\
&=&\left(
         \begin{array}{c}
        -0.55 \left(\arctan u^1_1 -\arctan u^2_1\right)+4.1 \left(u^1_1 -u^2_1\right)\\
        0.12\left(\sin u^1_2 - \sin u^2_2\right)+3.72 \left(u^1_2 - u^2_2\right)\\
         \end{array}
       \right)^\top \left(
         \begin{array}{c}
          u^1_1 - u^2_1 \\
           u^1_2 -u^2_2 \\
         \end{array}
       \right)\nonumber\\
&=& -0.55 \left(\arctan u^1_1 -\arctan u^2_1\right)\left(u^1_1 - u^2_1\right)+4.1 \left(u^1_1 - u^2_1\right)^2  +0.12 \left(\sin u^1_2 - \sin u^2_2\right)\left(u^1_2 - u^2_2\right) \nonumber\\
&&+3.72 \left(u^1_2 -u^2_2\right)^2 \nonumber\\
&=& -0.55 \frac{1}{1+\kappa_1^2} \cdot\left(u^1_1 - u^2_1\right)^2 + 4.1 \left(u^1_1 - u^2_1\right)^2 +0.12 \cos \kappa_2\cdot \left(u^1_2 - u^2_2\right)^2 +3.72 \left(u^1_2 -u^2_2\right)^2 \nonumber\\
&\geq& 3.55 \left(u^1_1 - u^2_1\right)^2 + 3.6 \left(u^1_2 -u^2_2\right)^2          \nonumber\\
&\geq& 3.55 \left(\left(u^1_1 - u^2_1\right)^2 +\left(u^1_2 -u^2_2\right)^2\right) \nonumber\\
&=& 3.55 \left\|u^1-u^2\right\|_2^2= m_{\widehat{S}} \left\|u^1-u^2\right\|_2^2,
\end{eqnarray*}
where $\kappa_i$ $(i=1,2)$ lies between $u^1_i$ and $u^2_i$. Therefore, $\widehat{S}$ is strongly monotone in its first variable with constant $m_{\widehat{S}}>o$, that is, assumption (A$'_6$) holds. And it follows immediately that
$$
 \widetilde{\lambda}=2\frac{L_F + L_g \eta_{\widehat{K}}+ \eta_g \frac{L_Q}{m_{\widehat{S}}}}{\Gamma(q+1)}T^q + L_{c_2}T+L_{c_1}T=0.9205 <1.
$$
Therefore, all assumptions of Theorem \ref{usc-PBFFDVI} hold, and so the solution mapping for PBFFDVI (\ref{ex2}) is u.s.c.
\end{example}

\section{Conclusions}
\noindent \setcounter{equation}{0}
In this work, we propose and investigate the BFFDVI, which provides a novel theoretical framework for characterizing fuzzy fractional boundary value problems constrained by VIs and supplies a new theoretical tool for studying boundary value problems of physical processes and materials with memory or anomalous phenomena in uncertain environments. We study the solvability and stability of the BFFDVI and present two numerical examples. The major contributions of this paper are outlined as follows: (i) the BFFDVI is proposed; (ii) the existence of solutions to the BFFDVI (\ref{BFFDVI}) is proved via the set-valued Krasnoselskii fixed point theorem; (iii) to analyze the stability of the BFFDVI (1.2), a new existence result for the parametrized BFFDVI (\ref{BFFDVI}), namely the PBFFDVI (\ref{PBFFDVI}), is established based on the set-valued contraction principle; (iv) the closedness of the solution mapping for the PBFFDVI (\ref{PBFFDVI}) is examined; (v) the u.s.c. of the solution mapping for the PBFFDVI (\ref{PBFFDVI}) is established; (vi) two numerical examples are provided to verify the validity of our main theoretical results. As noted previously, the semicontinuity of the solution mapping enables us to predict the variations in solutions of the BFFDVI (\ref{BFFDVI}) caused by parameter perturbations. Accordingly, exploring the l.s.c. of solution mapping for the PBFFDVI (\ref{PBFFDVI}) will be one of our directions for future research.

\end{document}